\documentclass[preprint,12pt]{elsarticle}

\usepackage{amsmath,amsfonts}
\usepackage{amsthm}
\usepackage{amssymb}
\usepackage{algorithmic}
\usepackage{array}
\usepackage{booktabs}
\usepackage{textcomp}
\usepackage{stfloats}
\usepackage{url}
\usepackage{verbatim}
\usepackage{graphicx}
\usepackage{pifont}
\usepackage{threeparttable}
\usepackage{tabularx}
\usepackage{soul}
\usepackage{xcolor}
\usepackage{longtable}
\usepackage{mathtools}
\usepackage{orcidlink}
\usepackage{hyperref}
\usepackage{cleveref} 
\usepackage{subcaption}  
\usepackage{graphicx}
\usepackage{enumitem}

\usepackage{fontawesome}
\usepackage{academicons}

\newcommand{\rv}[1]{\textcolor{black}{#1}}
\newtheorem{remark}{Remark}

\newtheorem{problem}{Problem}

\newtheorem{definition}{Definition}

\journal{Nuclear Physics B}

\begin{document}

\begin{frontmatter}


% \author{Name\corref{cor1}\fnref{label2}}
% \affiliation{organization={},

\title{Geometrically Explicit Cosserat-Rod Modeling with Piecewise Linear Strain for Complex Rod Systems}


\author{Lingxiao Xun, Brahim Tamadazte} %% Author name

\affiliation{organization={Sorbonne Université, CNRS UMR 7222, Inserm U1150, ISIR},%Department and Organization
            city={Paris},
            postcode={75005}, 
            country={France}}

%% Abstract
\begin{abstract}
This paper presents a geometrically explicit formulation for Cosserat rods that unifies configuration-space and strain-based representations within a single modeling framework. The proposed method uses nodal configurations on the Lie group \(SE(3)\) as generalized coordinates, while internal strains are reconstructed via a piecewise-linear parameterization. This hybrid design preserves the geometric rigor of Lie-group formulations and retains the locality, simplicity, and computational efficiency characteristic of strain-parameterized rod models. 

The formulation naturally avoids shear and membrane locking without additional stabilization techniques, and it accommodates arbitrary rod networks, closed-loop architectures, and gridshell-like structures through element-wise assembly. A Riemannian Newton solver is further developed to solve the equilibrium equations directly on \( SE(3) \), providing rapid convergence and consistent treatment of rotations. Numerical examples demonstrate that the method achieves high accuracy with only a few elements and generalizes seamlessly from single rods to complex multi-rod systems. These properties highlight the potential of the proposed formulation as a fast, robust, and scalable simulation tool for slender mechanisms and compliant structures.
\end{abstract}

% \includegraphics[width=\textwidth]{fig/abstract.png}

% \begin{highlights}
% \item A geometrically explicit Cosserat-rod formulation unifying strain and configuration representations.
% \item Internal strains reconstructed via piecewise linear parameterization without locking phenomena.
% \item Supports arbitrary rod networks, closed-loop topologies, and gridshell-like mechanisms.
% \item Riemannian Newton solver developed for efficient equilibrium computation on \( SE(3) \).
% \item Accurate and scalable performance demonstrated from single rods to complex compliant structures.
% \end{highlights}

%% Keywords

\end{frontmatter}



\section{INTRODUCTION}
% ---------
%
Soft and slender structures undergo large strains and rotations, making continuum beam or Cosserat-rod formulations a natural choice for their modeling~\cite{armanini2023soft, qin2024modeling}. Unlike rigid articulated mechanisms, these structures bend, twist, compress, and shear continuously along their length, often while interacting with complex environments. Accurate rod models are therefore essential for the design of soft slender structures, topology optimization, and high-fidelity simulation.

The finite element method (FEM) remains a versatile and powerful tool for continuum simulation and has been widely applied to soft materials. However, slender geometries combined with finite rotations and large strains typically require very fine spatial discretization to avoid spurious stiffness and numerical locking. This significantly increases computational cost and makes classical FEM less attractive for real-time or interactive applications. Beam-based formulations, such as Euler–Bernoulli and Timoshenko models~\cite{olson2020euler, lindenroth2016stiffness}, as well as geometrically exact Kirchhoff rod formulations~\cite{novelia2018discrete}, alleviate part of this burden by exploiting slender-body kinematics to reduce the number of degrees of freedom. Pseudo–rigid-body approaches~\cite{10.1115/1.3046148} further achieve reduced-order dynamics with only a few generalized coordinates, but they typically struggle to capture large spatial curvature, torsion–bending coupling, and complex loading conditions in strongly deformable structures.

Cosserat rod theory offers a rigorous continuum model with six strain and curvature components. Geometrically exact rod formulations from Simo~\cite{simo1985finite} account for $SO(3)$ rotations, where additive interpolation fails to preserve objectivity~\cite{crisfield1999objectivity}.  
Geodesic and helicoidal mappings~\cite{ sonneville2014geometrically}, together with modern Lie-group treatments~\cite{romero2018computing}, highlight the importance of geometric consistency under large rotations.

In recent years, a complementary research stream has adopted strain parameterized Cosserat models, treating strains as generalized coordinates to enable real-time and control-oriented formulations. This line includes piecewise constant (PCS) and piecewise linear (PLS) approximations~\cite{8500341, li2023piecewise}, geometric variable strain (GVS) models~\cite{9057619, boyer2020dynamics, mathew2025reduced}, and Newton–Euler strain integration~\cite{till2019real}. These methods excel in computational efficiency and analytic tractability. Nevertheless, strain integration from base to tip introduces challenges in spatial discretization, and extending these approaches to branched, interconnected, or closed chain rod systems \cite{armanini2021discrete} involves non-trivial constraint coupling, which presents additional difficulties for developing fully modular formulations suitable for general soft robotic, biological, or structural architectures.

In parallel, configuration space approaches, including isogeometric analysis (IGA) and Lie group finite elements, directly use pose as the generalized coordinates of the rod. IGA achieves high smoothness through NURBS directors~\cite{hughes2005isogeometric, weeger2016isogeometric}, supporting dynamics~\cite{weeger2018isogeometric}, contact~\cite{weeger2017isogeometric}, and fluid–structure interaction~\cite{agrawal2024efficient}. Nevertheless, these methods face challenges, such as the need for technical remedies to mitigate membrane and shear locking, and their dependence on spline continuity, which can complicate the construction of rod networks or hybrid soft structures, where changes in connectivity and topology are often required. 

\rv{Recent advances in geometrically exact beam discretizations have also emphasized high-order geometric constructions for manifold-valued interpolation, in particular for rotation fields. In this context, Greco and co-authors proposed spline- and Bézier-type interpolation schemes on the rotation manifold, together with De Casteljau-type algorithms, enabling high-order and geometrically consistent interpolation of orientations within nonlinear beam formulations~\cite{greco2024objective, greco2025spherical}.}

\rv{Alongside these manifold-based rotation constructions, a parallel line of work has developed geometrically exact beam/rod formulations using motion tensors and unit dual-quaternion representations to parameterize coupled translation--rotation kinematics. These frameworks are typically equipped with dedicated interpolation schemes and structure-preserving discretizations for large motions, and have been successfully applied to beams and rod networks~\cite{borri1994intrinsic,han2022configurational,han2019spectral,ghosh2008consistent}.}

Beyond $SE(3)$-based beam/rod formulations, closely related developments exist that rely on alternative but philosophically similar representations of rigid cross-section motion, notably motion-tensor descriptions and unit dual-quaternion parameterizations. These approaches also provide a unified treatment of translations and rotations and have been combined with structure-preserving discretizations and interpolation schemes for geometrically exact beams and rod networks

Although these research lines have matured, they typically emphasize complementary priorities: strain-driven models offer computational speed, compact representation, and real-time capability, whereas configuration-based models provide high-dimensional finite element simulations for complex structures. Nevertheless, seamlessly transitioning between single rods and intricate assemblies while preserving Lie group consistency, avoiding shear and membrane locking~\cite{mukherjee2001analysis}, and maintaining element-wise modularity remains a notable challenge. This motivates a formulation that combines the geometric fidelity of $SE(3)$ kinematics with the flexibility of strain reconstruction and the practicality of element-based assembly, enabling general rod networks and efficient simulation of large-deformation soft rod systems.
    
This work bridges the strain-based and configuration-space paradigms by introducing a geometrically explicit Cosserat rod formulation. We employ nodal poses on $SE(3)$ as primary variables and recover element-wise strain fields through a linear parameterization. This hybrid representation:
\begin{enumerate}
\item preserves $SE(3)$ invariance and naturally avoids shear and membrane locking without additional stabilization,
\item retains a clean element-wise structure, enabling arbitrary rod networks and closed loop geometries, and
\item accurately captures large rotations with only a small number of elements.
\end{enumerate}

By combining Lie group geometry with explicit spatial discretization and local strain reconstruction, the proposed approach unifies the key advantages of strain-parameterized and configuration-based models, offering a scalable framework for large-deformation slender structures and complex soft-rod architectures.

The rest of the paper is organized as follows.  Section~\ref{sec.cosserat} briefly reviews the classical Cosserat rod theory. Section~\ref{sec.element} presents the proposed discrete geometric explicit formulation, while Section~\ref{sec.statics} develops the static equilibrium conditions using Riemannian optimization. Section~\ref{sec.assembly} then describes the element residuals, tangent operators, and global assembly procedure, and Section~\ref{sec.ns} validates the method on single-rod benchmarks. Section~\ref{sec.na} demonstrates the capability of the approach on complex rod-network structures. Section~\ref{sec.conc} concludes the paper. 
%
% ------------
\section{COSSERAT ROD THEORY}\label{sec.cosserat}
% -------------
%
In this Section, we briefly review the fundamental theory of the Cosserat rod. 
% ----------
%
% --------
% \subsection{Continuous formulation}
% --------
%
In the Cosserat framework, the configuration of a deformable rod relative to the inertial frame is described by a position vector ${p}(s)\in\mathbb{R}^3$ and a rotational matrix ${R}(s)\in SO(3)$, parameterized by the abscissa material $s\in [0,L]$ along the rod, where \( s \) denotes the coordinate along the rod axis and \([0,L] \subset \mathbb{R}\) represents the reference domain of the length of the rod. Thus, the configuration of any cross-section of the soft rod can be defined as a curve ${g}(s): s\mapsto g(s)\in SE(3)$ with the homogeneous transformation matrix:
\begin{definition}[\textbf{\textit{Homogeneous transformation matrix}}]
The homogeneous transformation matrix for any cross-section along the soft rod can be defined as $$\forall s\in [0, L],\ {g}(s)=\begin{bmatrix}
	{R}(s)&{p}(s)\\
	{0}&1
\end{bmatrix}\in SE(3)$$
where ${p}(s)\in\mathbb{R}^3$ is the position vector and ${R}(s)\in SO(3)$ represents an orthonormal rotation matrix. $L$ is the total arc length of the soft rod.
\end{definition}

Subsequently, the strain and velocity are described in terms of the tangent space of the homogeneous transformation matrix.
\begin{definition}[\textbf{\textit{Strain of Cosserat rod}}]
	The strain in the body frame can be regarded as the left-trivialized tangent space of the homogeneous transformation matrix \textit{w.r.t.} space, i.e.,
	$$
		\hat{{\xi}}(s)={g}^{-1}{g}^{\prime}\in \mathfrak{se}(3)\simeq \mathbb{R}^6
	$$
\end{definition}
For simplicity, we use ${(\cdot)^\prime}$ to denote the partial derivative \textit{w.r.t.} space $\partial/\partial s$.

Based on the strain, the total internal energy of the rod can be expressed as:
\begin{equation}
\label{eq:total_energy}
U_\mathrm{int} =  \frac{1}{2}\!\int_0^{L} (\xi-\xi_0)^{\!\top}{K}(\xi-\xi_0)\,\mathrm ds.
\end{equation}
where $\xi_0$ represents the natural strain and the sectional stiffness matrix \( K \in \mathbb{R}^{6 \times 6} \) is defined as
\[
K = \mathrm{diag}(GJ_x,\, EJ_y,\, EJ_z,\, EA,\, GA,\, GA),
\]
where \(E\) is the Young’s modulus, \(G\) is the shear modulus, and \(A\) is the cross-sectional area. The quantities \(J_x\), \(J_y\), and \(J_z\) denote the second moments of area about the \(x\)-, \(y\)-, and \(z\)-axes, respectively.

Building on the preceding concepts and invoking D'Alembert's Principle, the equilibrium configuration of the rod corresponds to the stationary point of the total potential energy, i.e.,
\begin{equation}
\delta U_\mathrm{int} + \delta U_\mathrm{ext} = 0,
\end{equation}
with boundary conditions:
$
g(0) = g_0$ and (or) $g(L) = g_L
$, 
where $U_\mathrm{ext}$ represents the potential of external wrenches applied along the rod or at its boundaries.  

This continuous formulation represents a strain-based, energy-driven static model of the Cosserat rod. Various discretization strategies can be employed to solve the resulting system, and the choice of discrete degrees of freedom is crucial. Since the energy is directly expressed in terms of the strain field, a strain-based discretization is natural and facilitates the decoupling of axial and shear strains. On the other hand, boundary conditions are often defined in the configuration space, making configuration-based degrees of freedom more convenient for their enforcement and for modeling complex rod systems. In the following sections, we develop a framework that bridges these two approaches.    
This formulation preserves the full geometric structure of the Cosserat theory while enabling modular, element-wise computation, making it suitable for complex soft-rod architectures.
%
% ----------
\section{DISCRETE GEOMETRIC EXPLICIT FORMULATION}\label{sec.element}
% ----------
%
We now introduce the discrete geometric formulation used in this work. The continuous Cosserat rod is discretized into a sequence of nodes connected by finite segments (elements), with nodal poses in \(SE(3)\) as the main variables for representing the rod configuration in space.

Within each element, the internal strain field is integrated through a fourth-order Magnus expansion, yielding an explicit geometric expression that directly connects the nodal poses to the strain evolution along the rod. To mitigate shear locking in slender structures, a linearly varying strain distribution is assumed, preserving the natural separation between axial/shear strains and bending/torsional curvatures. At the global level, configuration continuity is enforced across elements, while strain continuity is relaxed, with linear variation maintained locally within each element.
%
% -------------------------------------------------------
\subsection{Linear Strain Element (LSE)}
% ---------------------------------------------------
%
Consider a single element of arc length $h$, parameterized by the coordinate $s\in[0,h]$.  
The spatial strain field within this element is assumed to vary linearly as
\begin{equation}\label{eq:linear_strain}
    \xi(s) = \bar{\xi} + (s - \tfrac{1}{2})\,\beta,
\end{equation}
where $\bar{\xi}\in\mathbb{R}^6$ represents the mean strain, and $\beta\in\mathbb{R}^6$ denotes the strain slope that captures the local variation of deformation along the element.  
%
% -------------------------------------------------------
\subsection{Fourth-order Magnus Approximation}
% ----------------------------------------------------
%
Integrating the kinematic relation $g'(s)=g(s)\hat{\xi}(s)$ under the linear strain field~\eqref{eq:linear_strain} leads to the exponential mapping between the nodal poses $g_a = g(0)$ and $g_b=g(h)$:
\begin{equation}
    g_b = g_a\,\exp\!\big(\Omega(h)\big),
\end{equation}
where $\Omega(h)\in\mathfrak{se}(3)$ is the integrated twist over the element via Magnus expansion~\cite{blanes2009magnus}.  

By applying the fourth-order Zanna–Magnus expansion and substituting~\eqref{eq:linear_strain}, one obtains a compact expression accurate up to $\mathcal{O}(h^5)$:
\begin{equation}\label{eq:magnus4}
    \Omega(h) = \big(h\,\mathbb{I} - \tfrac{h^{3}}{12}\operatorname{ad}_{\beta}\big)\bar{\xi}.
\end{equation}
where $\mathbb{I}$ denotes the identity matrix, and $\operatorname{ad}$ is the adjoint operator of the Lie algebra (see Appendix A)~\cite{8500341}. Appendix D presents the detailed derivation of eq.~\eqref{eq:magnus4}. Equation~\eqref{eq:magnus4} explicitly relates the integrated twist to the mean and slope of the strain distribution. 
Knowing the nodal poses $g_a$ and $g_b$ as well as the local slope $\beta$, the mean strain $\bar{\xi}$ can be retrieved without any iteration:
\begin{equation}\label{eq:explicit_meanstrain}
    \bar{\xi} = A^{-1}\,\mathrm{Log}\!\left(g_a^{-1}g_b\right),
    \qquad
    A = h\mathbb{I} - \tfrac{h^{3}}{12}\operatorname{ad}_{\beta}.
\end{equation}

The operator $\operatorname{ad}_{(\cdot)}$ denotes the adjoint action on $\mathfrak{se}(3)$, and $\mathrm{Log}(\cdot)$ is the matrix logarithm on $SE(3)$. 
This relation provides an explicit and geometrically consistent mapping between the nodal poses and the strain parameters $(\bar{\xi},\beta)$. The geometries of the element are fully determined by $\{g_a,g_b,\beta\}$.
%
% ----------------------------------------------------
\subsection{Kinematic Variations}
% ----------------------------------------------------
%
For static equilibrium and tangent stiffness evaluation, it is necessary to express the variation of the mean strain $\bar{\xi}$ in terms of the variations of the slope $\beta$ and the nodal pose right perturbations $(\delta\zeta_a,\delta\zeta_b)\in\mathfrak{se}(3)$ on the Lie group. 
The variation of ~\eqref{eq:magnus4} yields
\begin{equation}\label{eq:linearization}
    \delta\bar{\xi}
    = J_1\,\delta\zeta_a
    + J_2\,\delta\zeta_b+J_3\,\delta\beta,
\end{equation}
where the Jacobian matrices are
\begin{align}
    J_{1} &= A^{-1}\operatorname{dexp}^{-1}_{\Omega}\operatorname{Ad}_{\exp(-\hat\Omega)}, \label{eq:J1}\\[2pt]
    J_{2} &= A^{-1}\operatorname{dexp}^{-1}_{\Omega},\label{eq:J2}\\[2pt]
    J_{3} &= -\tfrac{1}{12}A^{-1}h^{3}\operatorname{ad}_{\bar{\xi}}. \label{eq:J3}
\end{align}

Here, $\operatorname{dexp}^{-1}_{\Omega}$ is the inverse differential of the exponential map, and $\operatorname{Ad}$ denotes the adjoint action of the Lie group (see Appendix)~\cite{8500341}. 
These Jacobians ensure the consistency of the linearization on the Lie group. They will be used to assemble the global residual and tangent matrices, which will be presented in the next section.

% ---------------------------------------------------
\subsection{Geometric Consistency and Modularity}
\begin{figure}[h]
	\centering
\includegraphics[width=0.75\textwidth]{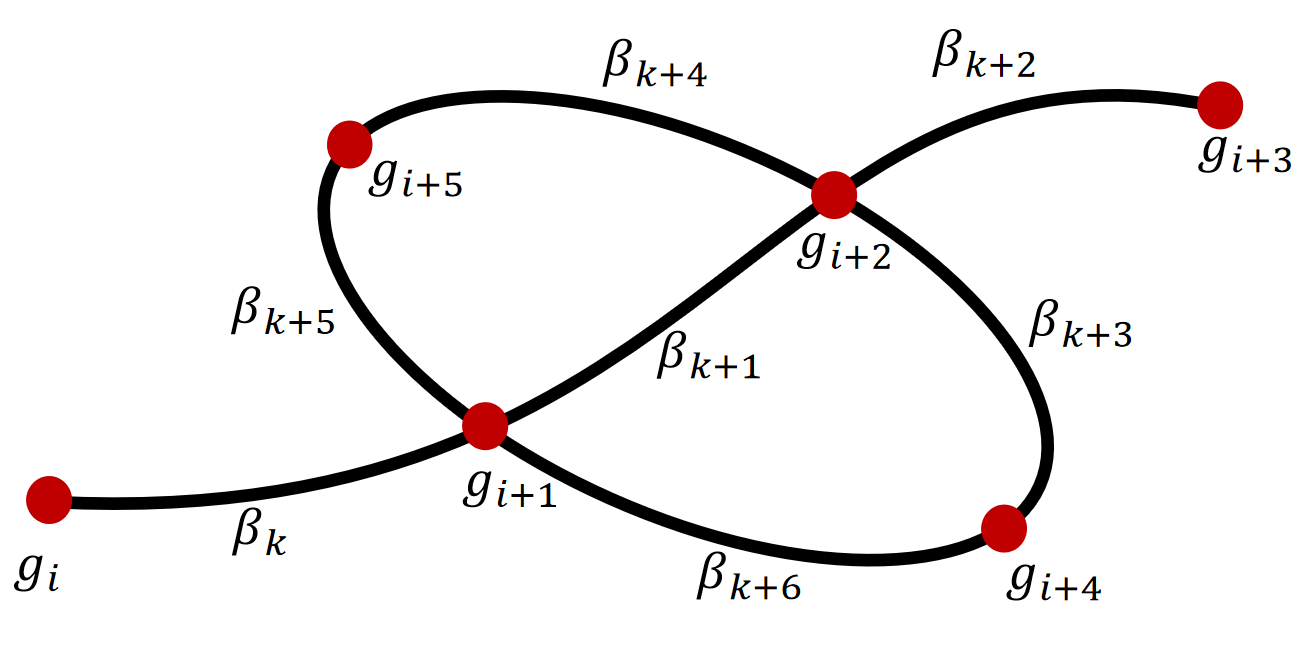}
\caption{ Topology of a rod network with $6$ nodes and $7$ elements.
The network consists of nodes and elements; each element connects two nodes and carries one independent strain-slope parameter. All incident elements share one node.}
	\label{fig:topoex}
\end{figure}

% ---------------------------------------------------
%
All quantities in the proposed formulation are expressed intrinsically on the Lie group $SE(3)$ and its algebra $\mathfrak{se}(3)$, ensuring exact treatment of rotations and translations. Each element is fully defined by its two nodal poses $(g_a,g_b)$ and one strain-slope parameter $\beta$, from which the internal strain field is reconstructed explicitly via~\eqref{eq:explicit_meanstrain}.  
Because the formulation does not impose inter-element strain continuity, the elements remain kinematically independent while maintaining geometric continuity of the nodal poses. This property enables the straightforward assembly of complex rod networks and seamless coupling with finite-element or multibody dynamics models, as illustrated in~\Cref{fig:topoex}. The resulting discrete system provides a compact and modular foundation for efficient computation, as detailed in the next section.
%
% -----------------
\section{STATIC EQUILIBRIUM via RIEMANNIAN OPTIMIZATION}\label{sec.statics}
% ----------------
%
In the discrete setting, if the rod system is divided into $N_e$ elements, its configuration is represented by the set of nodal poses
$\bar g=\{g_i\in SE(3)\mid i=0,\dots,N_n\}$
and by the set of strain-slope parameters
$\bar\beta=\{\beta_{e}\in\mathbb{R}^6\mid e=1,\dots,N_e\}$,
which together form the configuration vector
\[
q=(\bar g,\bar\beta)\in\mathcal M,
\qquad
\mathcal M = SE(3)^{N_n}\!\times\!\mathbb R^{6N_e}.
\]

The strain-slope parameters $\bar\beta$ encode the linear strain
distribution within each element, while the nodal poses $\bar g$
define the global geometry of the rod. The total potential energy of the discretized system can then be assembled as
\begin{equation}
U(q)
= \sum_{e=1}^{N_e}
  \Big(
    U_{\mathrm{int},e}(\bar\xi_e,\beta_{e})
    + U_{\mathrm{ext},e}(g_{a,e},g_{b,e})
  \Big),
\end{equation}
where $U_{\mathrm{int},e}$ and $U_{\mathrm{ext},e}$ respectively denote the internal and external potential energies of element $e$. Building on D'Alembert's Principle, the static equilibrium configuration corresponds to the stationary point of the total potential energy.

\begin{problem}[\textbf{\textit{Static equilibrium}}]
\label{pr.static}
Find the configuration set $q$ minimizing the total potential
energy:
\begin{equation}
    q
    = \operatorname*{arg\,min}_{q\in\mathcal M}
      U(q),
    \qquad
    \mathcal M = SE(3)^{N_n}\!\times\!\mathbb R^{6N_e}.
\end{equation}
\end{problem}

Note that since \( q \) lies on the product manifold \( \mathcal{M} \), the energy minimization is inherently a Riemannian optimization problem.
The following subsections introduce the essential mathematical concepts required for this formulation, including the retraction,
the corresponding directional derivative, and the definition of the Riemannian gradient~\cite{lang2012fundamentals}.
%
% ------------------------------------------------------------
\subsection{Basic Concepts of Riemannian Geometry}
\begin{definition}[\textbf{\textit{Retraction map}}]
Let $\mathcal G$ be a smooth manifold.
For any point $x\in\mathcal G$, let $T_x\mathcal G$ denotes the tangent space at $x$. A {retraction}
is a smooth mapping
\[
\mathrm{Ret}_x:T_x\mathcal G\to\mathcal G
\]
such that for any $v\in T_x\mathcal G$, the curve $c(t)=\mathrm{Ret}_x(tv)$ statisfies $c(0)=x$ and $\dot c(0)=v$.
% that satisfies $\mathrm{Ret}_x(0)=x$ and
% $\mathrm D\mathrm{Ret}_x(0_x)=\mathrm{id}_{T_x\mathcal G}$. 
\end{definition}

% The retraction provides a locally valid parameterization of $\mathcal G$
% around $x$.
% For instance, when $\mathcal G = SE(3)$, the retraction can be expressed through the right perturbation model
% \begin{equation}
% \mathrm{Ret}_x(v) = x\exp(\hat{\zeta}), 
% \qquad v = x\hat\zeta,\quad \zeta\in\mathfrak{se}(3).
% \end{equation}

The retraction maps tangent directions back to the manifold while preserving first-order geometry. For $SE(3)$, a common choice is the left trivialized retraction
\begin{equation}
\mathrm{Ret}_x(v) = x\,\exp(\widehat{\zeta}), 
\qquad 
v = x\hat\zeta,\quad \zeta\in\mathfrak{se}(3),
\end{equation}

Based on this mapping, let \( f:\mathcal{G} \to \mathbb{R} \) be a differentiable scalar field.
The {directional derivative} of \( f \) at \( x \) along any tangent vector \( v \in T_x\mathcal{G} \) 
is defined as:
\begin{equation}
\mathrm Df_x[v]
:= \left.\frac{\mathrm d}{\mathrm dt}\right|_{t=0}
    f\!\big(\mathrm{Ret}_x(t v)\big).
\label{eq:directional-derivative-general}
\end{equation}

This definition generalizes the classical Euclidean derivative
to curved manifolds and provides the basis for defining the following gradient.

\begin{definition}[\textbf{\textit{Riemannian gradient}}]
The {Riemannian gradient} of $f$ at $x$ is the unique tangent vector
$\operatorname{grad}f(x)\in T_x\mathcal G$ satisfying
\begin{equation}
\mathrm Df_x[v]
=\langle \operatorname{grad}f(x),\,v\rangle_x,
\qquad \forall v\in T_x\mathcal G,
\label{eq:riemannian-gradient-def}
\end{equation}
where $\langle\cdot,\cdot\rangle_x$ denotes the Riemannian metric on $\mathcal G$.
\end{definition}

When $\mathcal G = SE(3)$, we employ a left-invariant metric, which ensures that the inner product remains unchanged under left multiplication.
Accordingly, the {left-trivialized gradient} of $f$ at $g\in SE(3)$
is defined by
\[
\operatorname{grad}^L f(g):=
g^{-1}\operatorname{grad}f(g)
\;\in\;\mathfrak{se}(3),
\]
so that the directional derivative takes the compact form
\begin{equation}
\mathrm Df_g[g\hat\zeta]
=\zeta^\top\bigl(\operatorname{grad}^L f(g)\bigr)^\vee,
\qquad\forall\,\hat\zeta\in\mathfrak{se}(3).
\end{equation}
% where $\langle\cdot,\cdot\rangle_e$ denotes the standard inner product
% on the Lie algebra at the identity.
% 
% ----------
\subsection{Static Equilibrium Condition}
% ----------
%
The equilibrium configuration corresponds to a stationary point of
the potential energy $U(q)$, where its first variation vanishes.
By using the definition of the directional derivative on the product manifold
$\mathcal M = SE(3)^{N_n}\times\mathbb R^{6N_e}$, we can write
\begin{equation}
\delta U(q)
=\sum_{i=1}^{N_n} \mathrm D_{g_i}U[g_i\delta \hat\zeta_i]
 + \sum_{e=1}^{N_e} \mathrm D_{\beta_e}U[\delta\beta_e]
 = \delta q^\top r(q) = 0, 
\end{equation}
where $\delta q$ and $r(q)$ respectively denote the vectorized perturbation and the stacked form of the Riemannian gradient components:
\begin{equation}
\delta q=
\begin{bmatrix}
\delta \zeta_1\\[-2pt]\vdots\\[-2pt]\delta \zeta_{N_n}\\
\delta \beta_1\\[-2pt]\vdots\\[-2pt]\delta \beta_{N_e}
\end{bmatrix},
\qquad
r(q)=
\begin{bmatrix}
(\operatorname{grad}^L_{g_1}U)^{\!\vee}\\[-2pt]\vdots\\[-2pt](\operatorname{grad}^L_{g_{N_n}}U)^{\!\vee}\\
\nabla_{\beta_{1}}U\\[-2pt]\vdots\\[-2pt]\nabla_{\beta_{N_e}}U
\end{bmatrix}.
\label{eq:equilibrium-condition}
\end{equation}

Here $(\cdot)^{\vee}$ denotes the standard inverse of the hat operator
that maps $\mathfrak{se}(3)$ to $\mathbb R^6$.
The vector $r(q)$ therefore represents the static residual of the system, i.e., the collection of all nodal and strain-slope equilibrium errors.
The equilibrium configuration is characterized by the condition $r(q)=0$.
\begin{figure}[h]
	\centering
\includegraphics[width=0.75\textwidth]{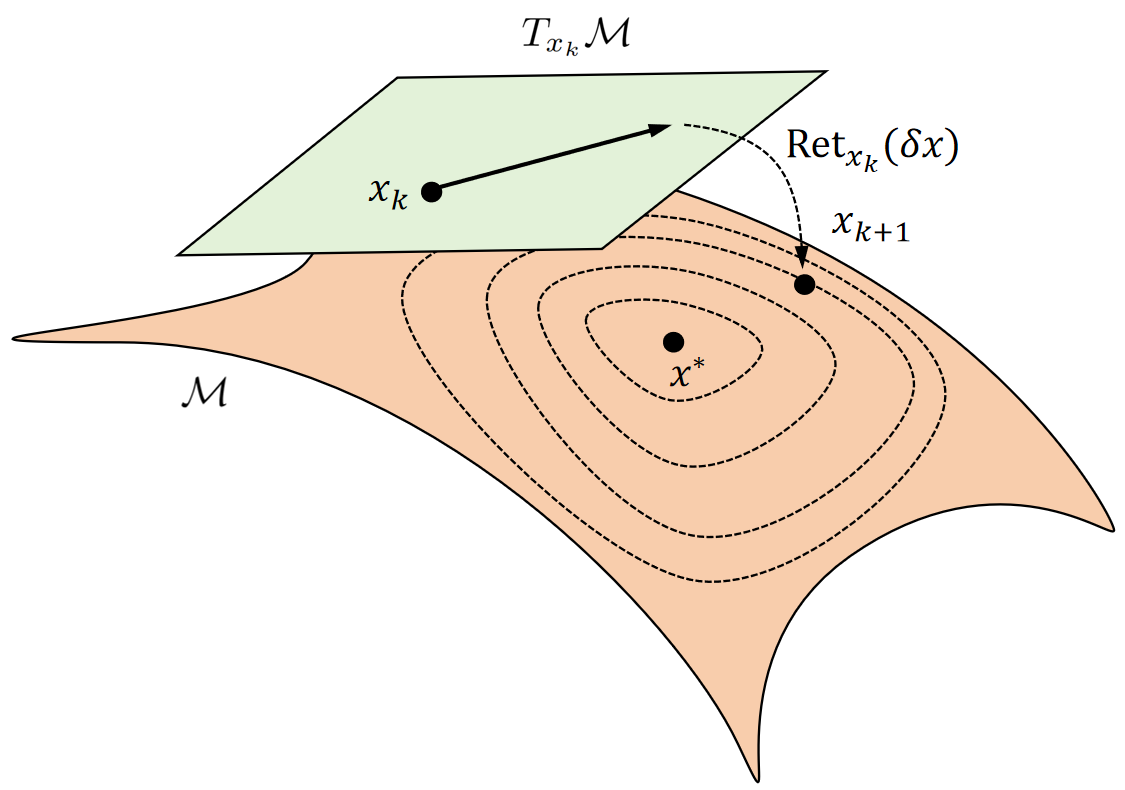}
\caption{Riemannian optimization scheme.
An update is computed in the tangent space \(T_{x_k}\mathcal{M}\) and mapped
back to the manifold via the retraction \(\mathrm{Ret}_{x_k}(\delta x)\), producing the next iterate \(x_{k+1}\) toward the solution \(x^\ast\).}
	\label{fig:ro}
\end{figure}
%
% --------
\subsection{Riemannian Optimization}
% --------
%
To numerically solve the equilibrium condition, a Riemannian Newton scheme~\cite{boumal2023introduction} is employed, as illustrated in~\Cref{fig:ro}. Denoting by $\delta q$ the incremental step at each iteration, the configuration of the Cosserat rod is updated as
\begin{equation}
g_i\leftarrow\mathrm{Ret}_{g_i}(\alpha g_i\delta\hat\zeta_i),
\qquad
\beta_e \leftarrow \beta_e+\alpha\,\delta\beta_e,
\end{equation}
where $\alpha>0$ denotes the step length.

The increment $\delta q$ is obtained from the linearized system
\begin{equation}
    \mathrm{Hess}\,U(q)\,\delta q=-r(q),
\end{equation}
where $\mathrm{Hess}\,U(q)$ denotes the Riemannian Hessian of $U(q)$ \cite{lang2012fundamentals},
defined as the covariant derivative of its gradient.
This iterative procedure converges to the equilibrium configuration
$q^\star$ for which the total potential energy is stationary.

The discrete equilibrium condition established in this section provides a macroscopic representation of the overall rod system. 
To render the formulation computationally tractable, we return to the element level, where the total potential energy \( U(q) \) is decomposed into element-wise contributions. We then derive explicit expressions for the element residuals and the tangent stiffness matrices (which serve as local approximations of the Hessian) and show how these quantities are assembled to construct the global residual vector and stiffness matrix of the entire rod system.
%
% ---------
\section{SYSTEM ASSEMBLY}\label{sec.assembly}
% ---------
%
\subsection{Elemental Energy and Equilibrium}\label{sec.eee}
% ----------------------------------------------
%
Substituting the linear strain field~\eqref{eq:linear_strain} into the energy integration~\eqref{eq:total_energy} yields the elastic energy stored in a single element:
\begin{equation}\label{eq:element_energy}
    U_{\mathrm{int},e} =
    \tfrac{1}{2}h(\bar{\xi}_e-\xi_{0,e})^{\!\top}{K}(\bar{\xi}_{e}-\xi_{0,e})
    + \tfrac{1}{24}h_e^{3}\,\beta_e^{\!\top}{K}\beta_e.
\end{equation}

The first term corresponds to the average strain deformation,
while the second term accounts for the strain-slop contribution.
For brevity, in the following~\cref{sec.eee} and \cref{sec.egn}, we assume that all quantities refer to the $e$-th element, and the subscript $e$ will therefore be omitted.

In statics, taking the first variation of~\eqref{eq:element_energy}
gives the elemental virtual work:
\begin{equation}\label{eq:element_virtualwork}
\begin{aligned}
    \delta U_{\mathrm{int}} &= 
    \mathrm{D}U_{\mathrm{int}_{g_{1}}}[\delta\zeta_a]
    + \mathrm{D}U_{\mathrm{int}_{g_{2}}}[\delta\zeta_b]
    + \mathrm{D}U_{\mathrm{int}_\beta}[\delta\beta]\\[2pt]
    &=
    \delta\zeta_a^{\!\top}R_1 +
    \delta\zeta_b^{\!\top}R_2 +
    \delta\beta^{\!\top}R_3,
\end{aligned}
\end{equation}
where $R_1$ and $R_2$ are the Riemannian gradients of the elemental energy with respect to $g_a$ and $g_b$, respectively, and $R_3$ is the gradient with respect to the strain-slope variable.
In physical terms, $R_1$ and $R_2$ correspond to the nodal residual wrenches acting at the left and right nodes, whereas $R_3$ represents the internal residual associated with the strain-slope parameter.
Substituting~\eqref{eq:linearization} into the variation of~\eqref{eq:element_energy} yields the explicit expressions:
\begin{align}
R_1 &= h\,J_1^{\!\top}\, K\,(\bar\xi-\xi_0),
\label{eq:R1}\\[2pt]
R_2 &= h\,J_2^{\!\top}\, K\,(\bar\xi-\xi_0),
\label{eq:R2}\\[2pt]
R_3 &= h\,J_3^{\!\top}\, K\,(\bar\xi-\xi_0)
      + \tfrac{1}{12}h^3\, K\,\beta.
\label{eq:Rk}
\end{align}
We denote the elemental residual of internal force as \begin{equation} 
F_{\mathrm{int}}=[R_1^{\top},R_2^{\top},R_3^{\top}]^{\top}\end{equation}.

\begin{remark}
(i) The nodal residuals~\eqref{eq:R1}–\eqref{eq:R2} represent the pullbacks of the mean-strain energy gradient along the Lie-group consistent directions defined by $J_1$ and $J_2$, thus ensuring kinematic consistency on $SE(3)$. (ii) The slope residual~\eqref{eq:Rk} consists of a coupling term $J_3^{\!\top}h K(\bar\xi-\xi_0)$
and a direct gradient-energy contribution $(h^3/12) K\beta$.  
\end{remark}

We assume that external loads act only at the nodal points. Let two external wrenches $F_a$ and $F_b$, expressed in the local frame, act on the nodes $g_a$ and $g_b$, respectively. The corresponding virtual work of the external loads is given by
\begin{equation}
\delta U_\mathrm{ext}
= -\delta\zeta_a^{\!\top}F_a
  -\delta\zeta_b^{\!\top}F_b.
\end{equation}

Hence, the elemental residual associated with the external forces can be written as
\begin{equation}
F_{\mathrm{ext}}
= [-F_a^{\!\top},\, -F_b^{\!\top},\, 0_{6\times1}]^{\!\top}.
\end{equation}

The total elemental residual is therefore obtained as
\begin{equation}
r_e = F_{\mathrm{int}} + F_{\mathrm{ext}}.
\end{equation}
%
% --------------------------------------------------
\subsection{Elemental Gauss--Newton Tangent} \label{sec.egn}
% ---------------------------------------------------
To solve the nonlinear equilibrium condition $r_e=0$ by a Newton-type iteration, the elemental residuals~\eqref{eq:R1}–\eqref{eq:Rk} must be linearized with respect to the incremental variables. The exact Riemannian Hessian of the elemental energy is generally cumbersome to compute, as it involves second-order derivatives of the Jacobians $J_{1,2,3}$ and exponential-map terms on $SE(3)$. To retain computational efficiency while preserving consistency, a Gauss–Newton approximation of the tangent operator is employed, with the Jacobians held constant at the current configuration.

Under this approximation, the differentials of the residuals read
\[
\begin{aligned}
\delta R_1 &= h\,J_1^{\!\top} K\,\delta\bar\xi, \\[2pt]
\delta R_2 &= h\,J_2^{\!\top} K\,\delta\bar\xi,\\[2pt]
\delta R_3 &= h\,J_3^{\!\top} K\,\delta\bar\xi + \tfrac{1}{12}h^3\, K\,\delta\beta.
\end{aligned}
\]

Collecting the coefficients of $(\delta \zeta_a, \delta \zeta_b, \delta \beta)$ 
yields the block tangent matrix 
$K^{\text{GN}} = \big(K^{\text{GN}}_{ij}\big)_{i,j=1}^{3}$, with entries
\begin{equation}\label{eq:GN_blocks}
K^{\text{GN}}_{ij} = 
\begin{cases}
h\, J_i^{\!\top} K J_j, & (i,j) \neq (3,3),\\[1.5mm]
h\, J_3^{\!\top} K J_3 + \tfrac{1}{12} h^3 \,K, & (i,j) = (3,3).
\end{cases}
\end{equation}

\begin{remark}
We note that the off-diagonal coupling blocks (e.g., $K_{13}$ and $K_{23}$ are not guaranteed to be full rank in general, since they stem from the pose--strain linearization and may weaken in special configurations. Nevertheless, the resulting Gauss--Newton tangent remains well-behaved in our setting because the diagonal blocks $K_{11}$, $K_{22}$, and, in particular, $K_{33}$ are coercive and well-conditioned, owing to the positive-definite Cosserat stiffness matrix and the explicit quadratic contribution of $\beta$ in the energy.
\end{remark}

% This Gauss–Newton tangent provides a symmetric,
% positive semi-definite approximation of the elemental Hessian
% and is used consistently in the numerical solver.
%
% ---------
\subsection{Gobal System} 
% ---------
The elemental residuals and Gauss--Newton tangents derived above 
serve as the building blocks of the global equilibrium system.
For each element $e$, let $A_e$ denote the Boolean assembly matrix
that maps the global perturbation vector $\delta q$ to the local variables $(\delta\zeta_{a,e},\delta\zeta_{b,e},\delta\beta_e)$. The matrix $A_e$ is constructed directly from the mesh topology,
encoding the nodal connectivity between adjacent elements.
The global residual vector $r(q)$ and global tangent matrix
$K^{\text{GN}}(q)$ are obtained by standard finite-element assembly:
\begin{equation}
r(q)=\sum_{e=1}^{N_e}A_e^{\!\top}r_e,
\qquad
K^{\text{GN}}(q)=\sum_{e=1}^{N_e}A_e^{\!\top}K_e^{\text{GN}}A_e.
\end{equation}

The constrained static equilibrium can thus be expressed as
\begin{equation}\label{eq:cons_tatic}
    r(q) = 0.
\end{equation}

Linearizing it yields:
\begin{equation}\label{eq:KKT}
K^{\text{GN}}\delta q = - r
\end{equation}
%
% ------------
\subsection{Global Newton Iteration}
% ------------
%
The assembled system defines the discrete equilibrium condition~\eqref{eq:cons_tatic}, which is solved using a Riemannian Newton iteration. At each iteration, the incremental update $\delta q$ is obtained by solving the linearized global system~\eqref{eq:KKT}, and the configuration variables are subsequently updated through the retraction maps:
\begin{equation}
g_i \leftarrow \mathrm{Ret}_{g_i}(\alpha g_i\delta\hat\zeta_i),
\qquad
\beta_e \leftarrow \beta_e + \alpha \delta\beta_e,
\end{equation}
where $\alpha$ denotes the step length.

The iteration proceeds until the residual norm $\|r(q)\|$
falls below a prescribed tolerance, yielding the equilibrium configuration of the complete rod system.
\begin{remark}\label{rmpcs}
The constant strain element (CSE), i.e, geodesic and helicoidal mappings~\cite{ sonneville2014geometrically}, can be regarded as a special case of the linear strain element, where the strain slope $\beta$ is set to zero. In the model, this simply corresponds to removing all terms associated with $\beta$.
\end{remark}
% 
% ----------------
\section{MODEL VALIDATION}\label{sec.ns}
% ----------------
In this section, we present a single-rod simulation to validate the performance of the proposed modeling method. All simulations were conducted within the MATLAB environment on a CPU {13th Gen Intel\textsuperscript{\textregistered} Core\textsuperscript{TM} i7-13850HX @ 2.10\,GHz}.
%
% -------------
% \subsection{Accuracy Validation}\label{sec.ava}
% ------------
%
This section validates the accuracy of the proposed geometric discretization through comparisons against reference solutions on two representative benchmarks: a planar cantilever under tip forces and a genuinely three-dimensional case involving bending--torsion coupling. 
% ------
\subsection{In-plane Validation}\label{sec.ava}
\begin{figure}[!h]
	\centering
\includegraphics[width=0.99\textwidth]{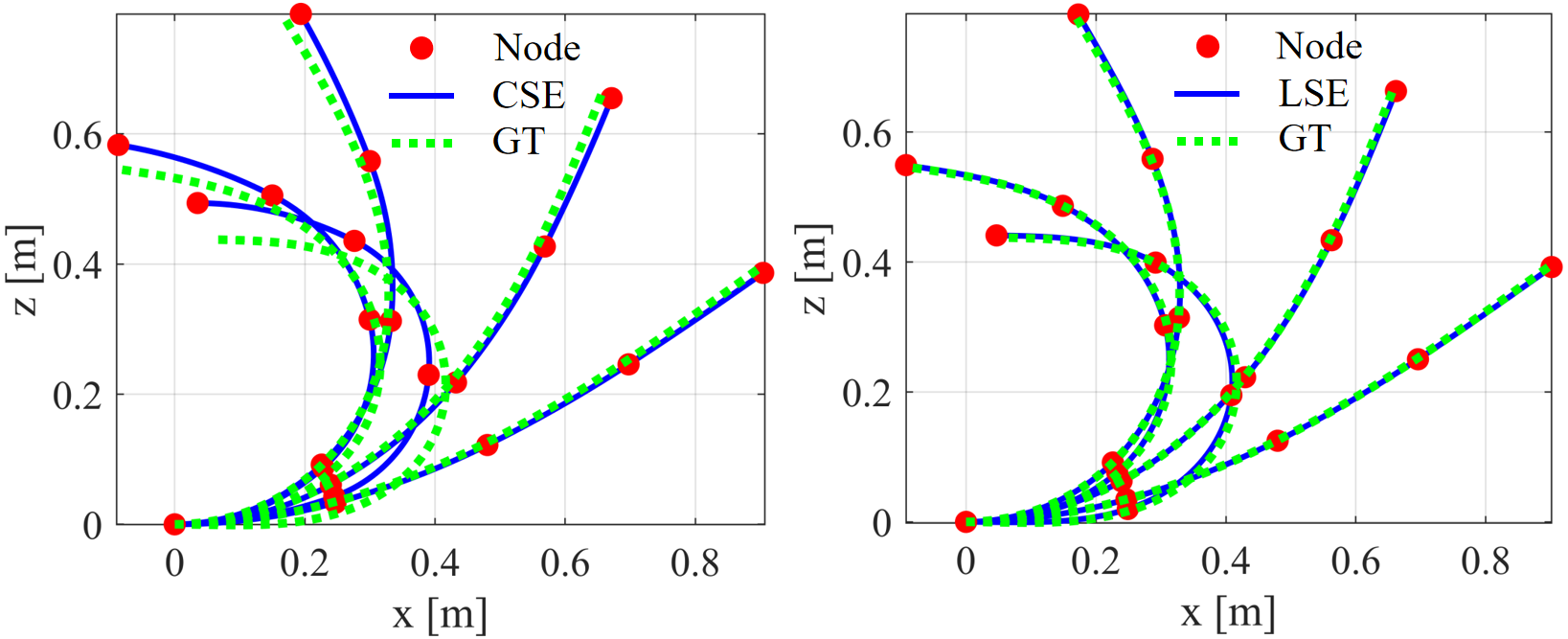}
\caption{Comparison of deformation results between the proposed two elements (CSE and LSE) and the ground truth (GT). }
	\label{fig:comparepls}
\end{figure}

This test evaluates the accuracy of the proposed modeling method by comparing simulation outcomes with ground-truth data. The experimental setup, illustrated in~\Cref{fig:comparepls}, involves a soft rod positioned horizontally with one end fixed at the origin. A force of 0.25~N, 0.5~N, 1~N, 2~N, and 3~N is applied along the rod's body-frame z-axis at its right end. The physical parameters of the soft rod used in the simulation are detailed in~\Cref{tab:cosserat_params}.
\begin{table}[!h]
  \centering
  \small
  \caption{Physical Parameters}
  \label{tab:cosserat_params}
  \begin{tabular}{l@{\hskip 100pt}c}
    \toprule
    \textbf{Parameter}        & \textbf{Value / Unit}            \\
    \midrule
    Length                   & $1$ m                      \\
    Bending stiffness          & $0.2$ Nm$^2$     \\
    Shear and axial stiffness          & $1000$ N     \\
    \bottomrule
  \end{tabular}
\end{table}

In the simulation, the rod was discretized into $4$ elements using both the constant strain element and linear strain element formulations. The constant strain element is a special case of a linear strain element, as detailed in~\cref{rmpcs}. The ground-truth solution is obtained by solving the static Poincaré equations of the Cosserat rod with the shooting method. 

\Cref{fig:comparepls} shows the simulation results. The proposed method exhibits close agreement with the ground truth, maintaining a relative error below \(1\%\). No locking behavior is observed, even with a coarse discretization of only \(4\) elements.

\rv{
\begin{remark}
In this work, the rest geometry is represented via the intrinsic strain $\xi_{0}$ appearing in the internal energy through $(\xi-\xi_{0})$. For initially curved rods, we assume $\xi_{0}$ to be piecewise constant per element and compute $\xi_{0,e}$ directly from the prescribed initial nodal poses $\{g_i^{0}\}$ using the same pose--strain mapping as in the proposed formulation. This construction applies to a single rod. For multi-rod systems with crossing nodes (e.g., gridshells), a dedicated treatment of nodal poses at intersections is required and will be addressed in future work.
\end{remark}}
%
% -----------
\subsection{Out-of-plane Validation}\label{sec.ava}
\begin{figure}[!t]
	\centering
	\includegraphics[width=0.5\textwidth]{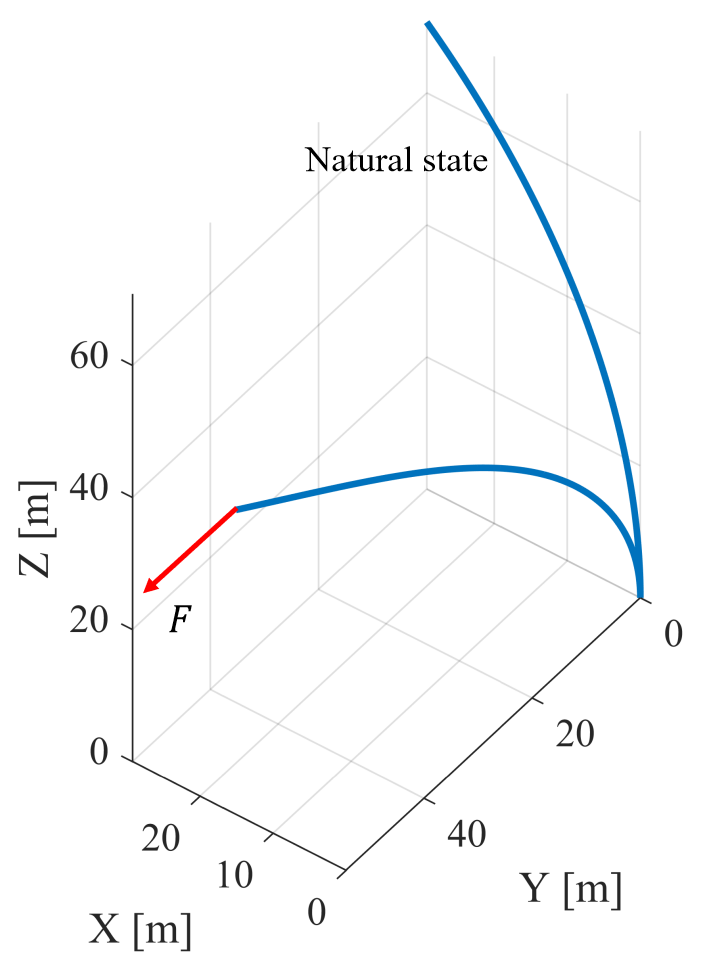}
	\caption{Cantilever 45-degree bend subjected to a fixed load.}
	\label{fig:outplane}
\end{figure}

\rv{The second benchmark reproduces the standard three-dimensional test case proposed in \cite{sonneville2014geometrically,meier2018geometrically}. The rod is initially curved as a $45^\circ$ bend with radius of curvature $100~\mathrm{m}$, length $L=25\pi~\mathrm{m}$, and circular cross-section of radius $r=1~\mathrm{m}$. The material is linear elastic with Young's modulus $E=10^{6}~\mathrm{MPa}$ and Poisson's ratio $\nu=0$. A tip force of magnitude $F=600~\mathrm{N}$ is applied along the global $y$-axis. Because the deformed configuration is non-planar, this benchmark activates bending--torsion coupling and provides a more stringent validation than purely planar tests. The natural and deformed shapes are shown in~\Cref{fig:outplane}.}

\subsubsection{Accuracy versus Degrees of Freedom}\label{sec.sss}
\rv{The proposed LSE is compared with the constant-strain element (CSE) used in \cite{sonneville2014geometrically} and the classical three-node Simo--Reissner (SR) beam element reported in \cite{meier2018geometrically}. To ensure a fair comparison, the total number of DoFs is matched between CSE, LSE, and SR, as summarized in Table~\ref{tab:dofs}. A dense-mesh numerical solution ($1000$ elements) is used as reference. }
\begin{table}[!h]
\centering
\renewcommand{\arraystretch}{1.3}
\caption{Meshes and degrees of freedom (DoFs) for the 45-degree bend benchmark, with matched DoFs between CSE and LSE.}
\label{tab:dofs}
\begin{tabular}{c c c c c c}
\hline
\hline
 & & Test 1 & Test 2 & Test 3 & Test 4\\
\hline
CSE & Elements & $4$  & $8$  & $12$ & $16$ \\
    & DoFs     & $30$ & $54$ & $78$ & $102$ \\
\hline
LSE & Elements & $2$  & $4$  & $6$  & $8$ \\
    & DoFs     & $30$ & $54$ & $78$ & $102$ \\
    \hline
SR & Elements  & $2$  & $4$  & $6$  & $8$ \\
    & DoFs     & $33$ & $57$ & $81$ & $105$ \\
\hline
\hline
\end{tabular}
\end{table}

\begin{figure}[!t]
	\centering
	\includegraphics[width=0.9\textwidth]{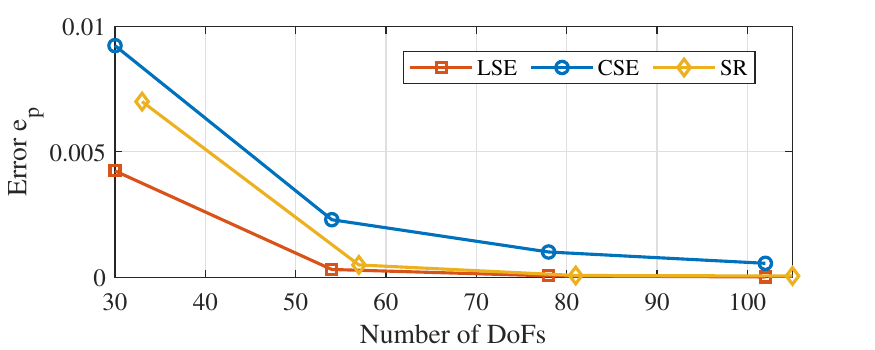}
	\caption{Tip-displacement error $e_p$ versus the total DoFs for CSE and LSE (reference: dense mesh with $1000$ elements).}
	\label{fig:outplanedoferror}
\end{figure}

\rv{The displacement error is defined as}
\rv{\[
e_p \;=\; \frac{1}{u_{max}}\sqrt{\frac{1}{L}\int_0^L\big\|p(s)-p_r(s)\big\|^2},
\]}
\rv{where $p(s)$ denotes the computed position and $p_r(s)$ denotes the reference position. $u_{max}$ denotes the maximal displacement. \Cref{fig:outplanedoferror} reports $e_p$ as a function of the total DoFs for both discretizations.}

\rv{The results show that LSE consistently yields smaller errors than CSE for the same DoF budget. As the DoFs increase, both methods exhibit monotonic decay in the displacement error, with LSE decreasing more rapidly, indicating improved accuracy-per-DoF efficiency. This improvement stems from LSE's ability to represent intra-element strain variations, thereby reducing the approximation error associated with piecewise-constant strain assumptions.}
\subsubsection{Path Independence}
\rv{To verify the path independence of our proposed method, we consider the same final load level $F=600~\mathrm{N}$ and solve the problem using three different loading paths:
(i) a single-step application ($F=600~\mathrm{N}$ applied directly),
(ii) a 10-step linear ramp $F_k = \frac{k}{10}F$ ($k=1,\dots,10$),
and (iii) a 10-step smooth ramp based on a sine function, $F_k = \sin\!\big(\frac{k\pi}{20}\big)\,F$ ($k=1,\dots,10$), which also reaches $F$ at the last step.
For each path, the Newton solver is run to the same convergence tolerance at every load increment.}

\rv{We then compare the final converged tip position with that obtained from single-step loading. The differences for the 10-step linear ramp and the 10-step sine ramp are $3.52\times 10^{-15}$ and $5.71\times 10^{-15}$, respectively, i.e., essentially at machine precision. This confirms that the computed solution is path-independent for this static benchmark and provides additional evidence of the consistency of the proposed numerical formulation.}

\subsubsection{Convergence Rate}\label{sec.ssr}
\rv{To quantify convergence rate with a standard and objective metric, an energy-type error is defined from the strain-field discrepancy,}
\rv{\[
e_E \;=\; \int_0^L\big(\xi(s)-\xi_r(s)\big)^\top K_\xi \,\big(\xi(s)-\xi_r(s)\big)\,\mathrm{d}s,
\]}
\rv{where $K_\xi$ denotes the Cosserat stiffness matrix in the body frame, $\xi(s)$ is the strain predicted by the discretization under test, and $\xi_r(s)$ is a reference strain field obtained from a mesh-converged high-accuracy solution ($1000$ elements). This stiffness-weighted metric corresponds to a strain-energy norm and provides a meaningful basis for assessing convergence rates.}

\begin{figure}[!h]
	\centering
	\includegraphics[width=0.8\textwidth]{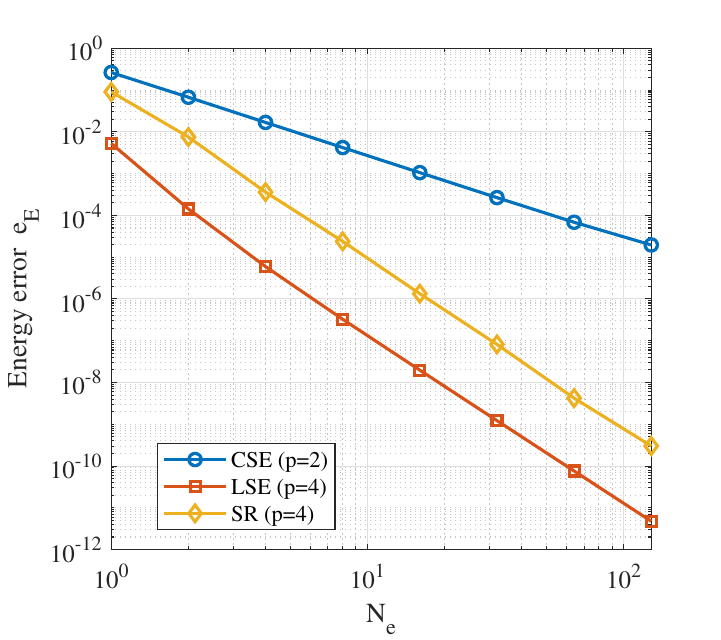}
	\caption{Log--log plot of the stiffness-weighted strain-field error $e_E$ versus the number of elements $N_e$, comparing CSE, LSE, and the three-node Simo--Reissner element.}
	\label{fig:convergence_rate2}
\end{figure}

\rv{A systematic mesh-refinement study is performed for CSE, LSE, and SR with $N_e = 1,\,2,\,2^2,\,\ldots,\,2^7$ elements.  \Cref{fig:convergence_rate2} shows the log--log plot of $e_E$ versus the number of elements $N_e$, from which convergence rates can be directly assessed. The observed trends are consistent with the respective assumed-field interpolations and confirm that LSE achieves a faster reduction (convergence order $p=4$) of the strain-field error under refinement compared with CSE (convergence order $p=2$) for this benchmark.}

%
% ------------
\subsection{Convergence Behavior}\label{cbs}
% ------------
%
Here, we evaluate the convergence behavior of the proposed static solver. The test case uses the same soft rod from the previous section, with one end fixed, and a force of 1 N applied at the right end of the rod along the body-frame z-axis. 

\begin{figure}[!h]
	\centering
\includegraphics[width=0.9\textwidth]{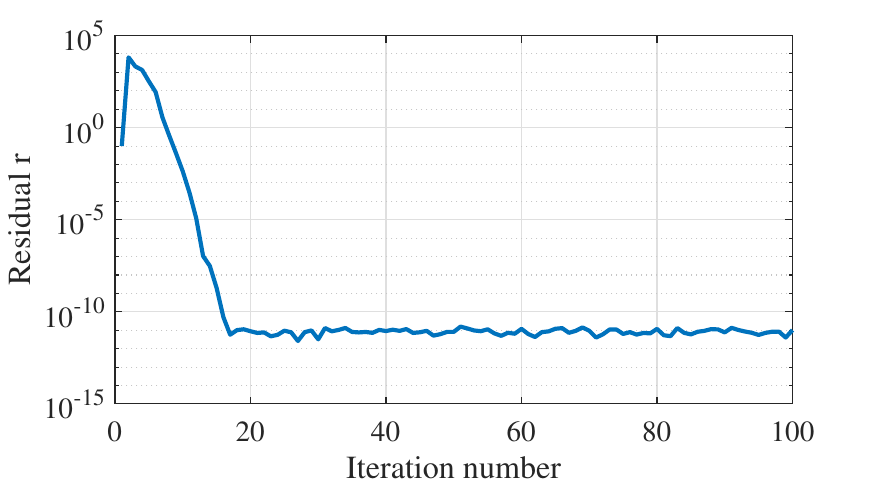}
	\caption{Convergence behavior of the static solver. }
	\label{fig:iterations}
\end{figure}

The initial guess for the shape of the soft rod is its naturally straight configuration. \Cref{fig:iterations} illustrates the evolution of the residual norm $\|{r}\|$ over successive iterations of the static solver. The solver exhibits rapid convergence, with the residual dropping below $10^{-3}$ within approximately ten iterations. Subsequently, the residual continues to decrease and stabilizes at around $10^{-11}$, indicating the robustness and numerical stability of the proposed method. \rv{The reason the convergence no longer decreases is that we employ a Gauss–Newton approximation of the Hessian in the optimization. As the iterates approach the solution, the convergence rate may deteriorate from second-order to first-order, or even stall altogether. One possible remedy is to derive and use the exact analytical Hessian, which we plan to address in future work.}

\rv{\subsection{Locking-free Benchmarks}}
\label{subsec:lockingfree}

\rv{In slender rod simulations, certain discretizations may exhibit artificial stiffening in bending-dominated regimes. In particular, an undesirable behavior is the appearance of {spurious axial (tensile) strain} caused by an artificial coupling between bending and axial deformation; similarly, in the thin-beam limit, {shear locking} may lead to an overly stiff response under transverse loading. This subsection reports two standard numerical benchmarks designed to verify that the proposed geometric discretization remains free of such locking artifacts in the regimes considered.}

\begin{figure}[!h]
	\centering
	\includegraphics[width=0.75\textwidth]{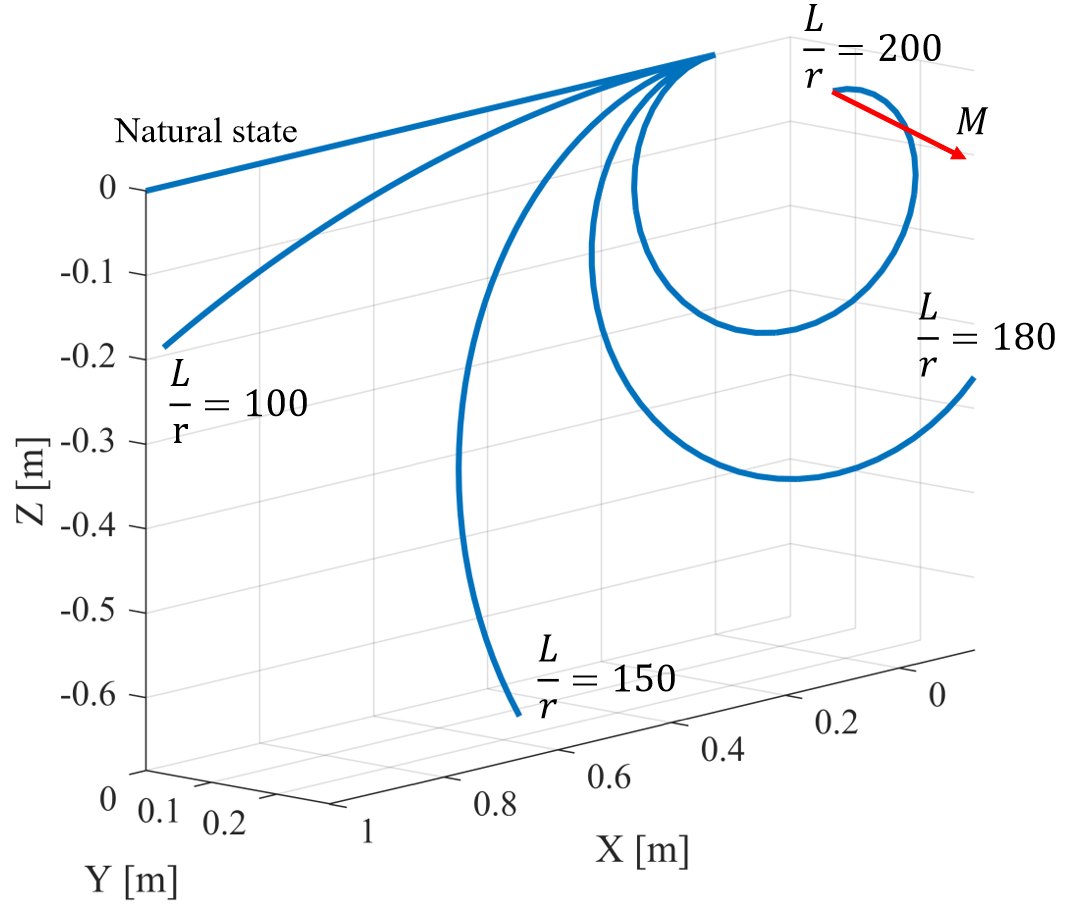}
	\caption{{Helical deformations of the rod under a constant end moment applied at the free tip, shown for different slenderness ratios $L/r=100,150,180,200$.}}
	\label{fig:bendingall}
\end{figure}

\rv{\subsubsection{Pure-bending Patch Test}}
\rv{A pure-bending patch test is first considered. The rod is clamped at one end in the reference configuration, and a constant end moment (expressed in the global frame) is applied at the free tip,
$
M = [\,5,\ 20,\ 0\,]~\mathrm{mN\,m}.
$
The rod is discretized with $4$ elements. Several slenderness ratios are investigated by varying the radius $r$, namely $L/r\in\{100,150,180,200\}$. The resulting deformed configurations are shown in \Cref{fig:bendingall}. The response is bending-dominated and does not display spurious axial effects in this regime. In addition, the recovered bending behavior is consistent with the analytical helical solution, supporting the correctness of the formulation in pure bending.}

\begin{figure}[!h]
	\centering
	\includegraphics[width=0.9\textwidth]{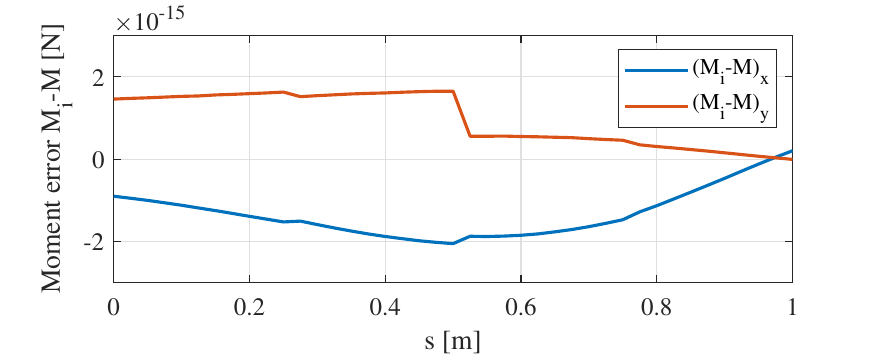}
	\caption{Error in the internal bending-moment distribution (global frame) with respect to the analytical solution for the pure-bending test ($L/r=200$, $4$ elements).}
	\label{fig:bendingmoment}
\end{figure}

\rv{To quantify the accuracy, the internal bending-moment distribution predicted by the formulation is compared with its theoretical value. For a rod subjected to a pure end moment, the internal moment $M_i$ is theoretically constant along the arc length and equal to the applied external moment $M$ at every cross-section. \Cref{fig:bendingmoment} reports the difference $M_i-M$ along the arc length for the case $L/r=200$ discretized with four elements (moments expressed in the global frame). The error is close to machine precision, indicating that the method achieves high accuracy even with a very coarse mesh in pure bending conditions.}

\begin{figure}[!t]
	\centering
	\includegraphics[width=0.9\textwidth]{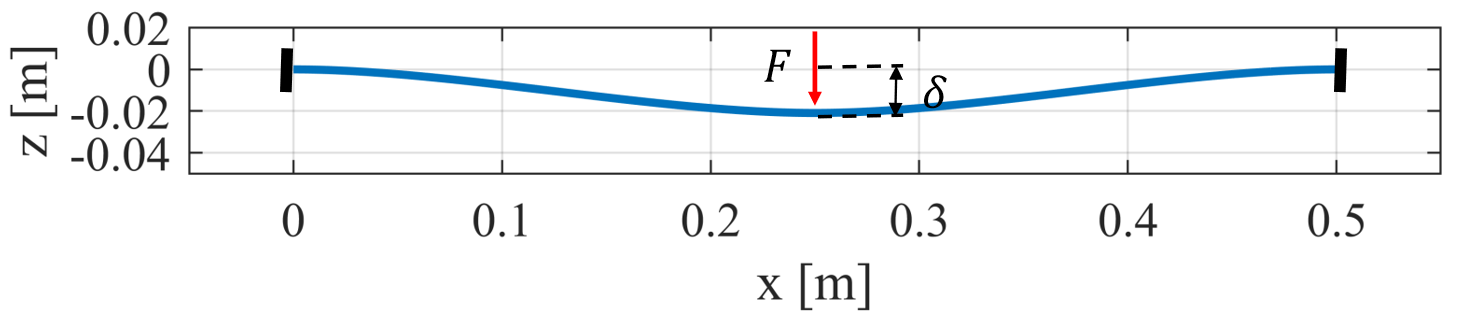}
	\caption{{Schematic of the clamped--clamped beam under a mid-span concentrated load.}}
	\label{fig:shearlocking2}
\end{figure}

\rv{\subsubsection{Clamped--clamped Beam under a Mid-span Concentrated Load}}
\rv{A classical clamped--clamped beam benchmark under a mid-span concentrated load is then considered to assess robustness in the thin-beam limit and to detect possible shear-locking induced stiffening. The rod is fixed at both ends, and a transverse load is applied at the mid-span along the $z$ direction, as illustrated in \Cref{fig:shearlocking2}. Three slenderness ratios are tested by varying the radius, namely $L/r = 50, 100,$ and $200$. For each case, several mesh resolutions ($N_e=8, 16, 32$) are considered, and a fine-mesh-converged numerical solution ($N_e=200$) is used as the reference.}

\rv{Table~\ref{tab:benchmark2_conv} reports the mid-span displacement $\delta$ and the corresponding relative error with respect to the reference solution. A rapid reduction in displacement error is observed with mesh refinement across all slenderness ratios. For the most slender case ($L/r=200$), the relative error decreases from $9.78\%$ ($N_e=8$) to $0.13\%$ ($N_e=16$), and further to $2.5\times 10^{-2}\%$ ($N_e=32$). These results indicate that the proposed method remains accurate in the thin-beam regime considered here and does not exhibit shear-locking-induced stiffening. For a fixed mesh resolution, the error may increase with slenderness because deformation becomes more localized, requiring a finer discretization to achieve the same accuracy.}

\begin{figure}[!t]
	\centering
	\includegraphics[width=0.9\textwidth]{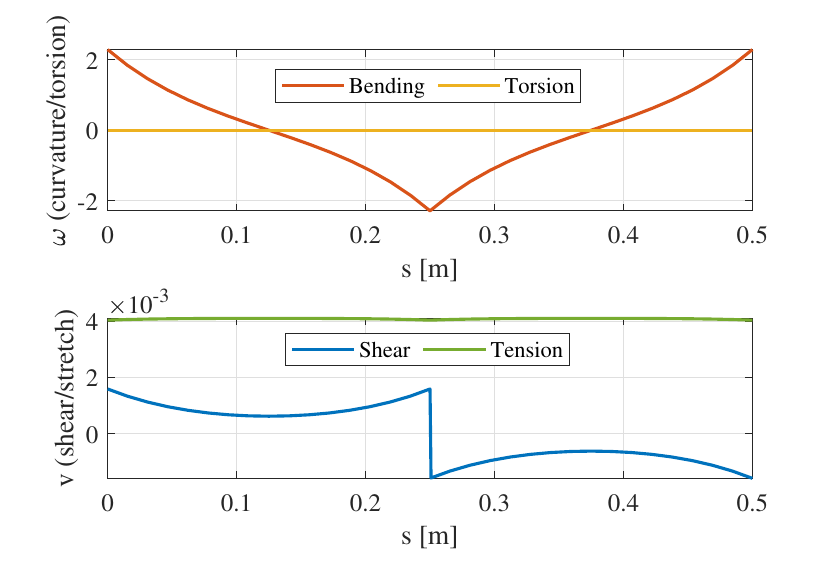}
	\caption{Distribution of strain in the body frame ($L/r=50$, $N_e=32$).}
	\label{fig:shearlockingstrain}
\end{figure}

\begin{figure}[!t]
	\centering
	\includegraphics[width=0.9\textwidth]{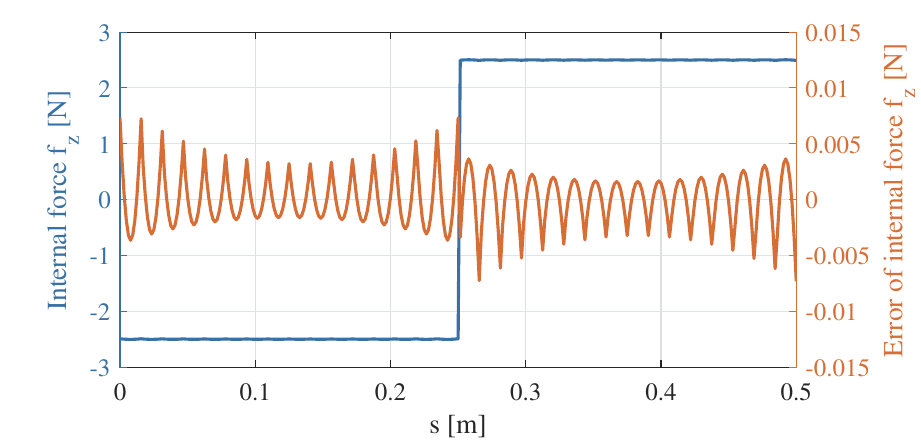}
	\caption{Distribution of internal force $f_z$ and its error compared with the analytical value in the global frame ($L/r=50$, $N_e=32$).}
	\label{fig:shearlockingforce}
\end{figure}

\rv{Beyond global quantities, internal-field distributions are also examined to rule out spurious oscillations under refinement. \Cref{fig:shearlockingstrain} shows representative strain components in the body frame for the case $L/r=50$ with $N_e=32$, exhibiting a smooth spatial profile without oscillations. \Cref{fig:shearlockingforce} reports the corresponding internal force component $f_z$ in the global frame; it agrees well with the analytical distribution, with an absolute discrepancy within $0.01~\mathrm{N}$ in the present settings, and the distribution remains smooth along the arc length.}

\begin{table}[!t]
\centering
\renewcommand{\arraystretch}{1}
\caption{Clamped--clamped beam under a mid-span concentrated load.
For each radius $r$, several mesh resolutions are considered to assess convergence. The reference solutions are obtained from a fine mesh-converged numerical solution ($200$ elements).}
\label{tab:benchmark2_conv}
\begin{tabular}{c c c c c}
\hline
\hline
Radius $r$ &
$L/r$ &
$N_e$ &
Displacement $\delta$ &
Relative error (\%) \\
\hline
$1\times10^{-2}$ & $50$ & $8$  & $2.10\times10^{-2}$ m & $0.31$ \\
                 &      & $16$ & $2.09\times10^{-2}$ m & $1.9\times10^{-2}$ \\
                 &      & $32$ & $2.09\times10^{-2}$ m & $1.1\times10^{-3}$ \\
\hline
$5\times10^{-3}$ & $100$ & $8$  & $2.17\times10^{-2}$ m & $1.28$ \\
                 &       & $16$ & $2.14\times10^{-2}$ m & $8.2\times10^{-2}$ \\
                 &       & $32$ & $2.14\times10^{-2}$ m & $1.7\times10^{-2}$ \\
\hline
$2.5\times10^{-3}$ & $200$ & $8$  & $2.31\times10^{-2}$ m & $9.78$ \\
                   &       & $16$ & $2.10\times10^{-2}$ m & $0.13$ \\
                   &       & $32$ & $2.10\times10^{-2}$ m & $2.5\times10^{-2}$ \\
\hline
\hline
\end{tabular}
\end{table}

%
% ---------
\section{NUMERICAL APPLICATIONS}\label{sec.na}
% ---------
This section provides numerical examples involving rod networks with various geometries. The aim is to evaluate the formulation on structures with multiple Cosserat rods intersecting, including setups relevant to soft-rod systems with lattice reinforcements or flexible shell-like components. We begin with planar and spatial rod grids, followed by a hemispherical lattice representing a discrete shell. All simulations are conducted in MATLAB using a unified $SE(3)$-based formulation and solver.
%
% -------
\subsection{Rod Networks}
% -------
Rod-based network structures frequently arise in soft structures, including deployable continuum truss mechanisms, and architected bio-inspired structures that emulate shell-like compliance via slender elastic members. These systems rely on large deformations, looped connectivity, and distributed strain-energy transfer, which pose challenges for strain-integrated Cosserat formulations that must enforce compatibility constraints.
\begin{figure}[!h]
    \centering
    \includegraphics[width=0.8\textwidth]{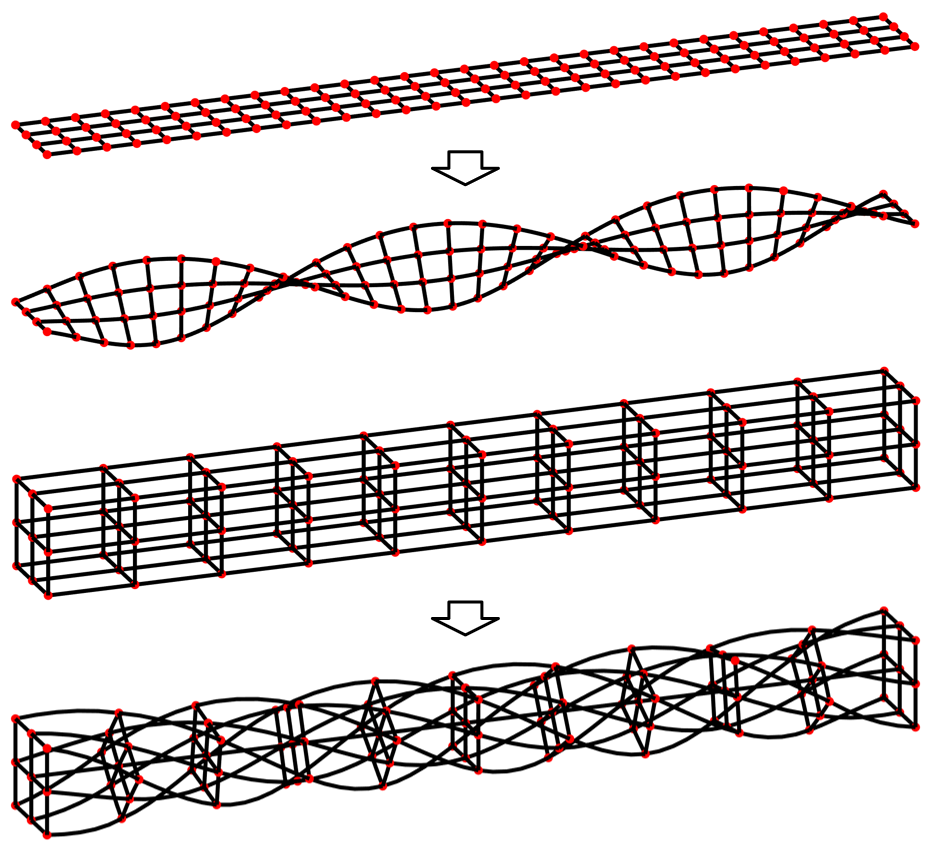}
    \caption{
    Deformed configurations of a planar lattice and a spatial rod network, modeled as Cosserat-rod networks
    }
    \label{fig:plat}
\end{figure}

\rv{The proposed \(SE(3)\)-based formulation uses nodal poses as the primary degrees of freedom, so rod networks, including closed loops, can be modeled without introducing additional compatibility or loop-closure constraints. Geometric continuity is enforced by construction through shared nodal poses, while internal strains are reconstructed locally at the element level. This structure enables efficient simulation of closed-loop lattices and gridshell-like rod assemblies encountered in soft continuum robots and mechanically intelligent metamaterials.}

We first examine a planar lattice composed of \(120\) nodes and \(206\) linear-strain rod elements. The left boundary is clamped, and the right boundary undergoes a prescribed twist around the longitudinal axis, producing a globally coupled torsional mode within the network (\Cref{fig:plat}).  Rod mechanical parameters match those used in~\cref{sec.ava}.

We then consider a 3D truss composed of \(99\) nodes and \(220\) rod elements. The left face is fixed, and the right face is twisted, yielding smooth, distributed deformation across a spatial structure with multiple closed loops. 

Both examples remain numerically stable and convergent under large rotations and kinematic coupling, highlighting the robustness of the Riemannian Newton solver for soft lattice systems.
%
%
% --------
\subsection{Parallel Mechanisms}
% ---------------
%
To further evaluate the efficiency, generality, and robustness of the proposed modeling and simulation framework, we simulate a representative and mechanically complex parallel soft-rigid structure. 
\begin{figure}[h]
	\centering
\includegraphics[width=0.9\textwidth]{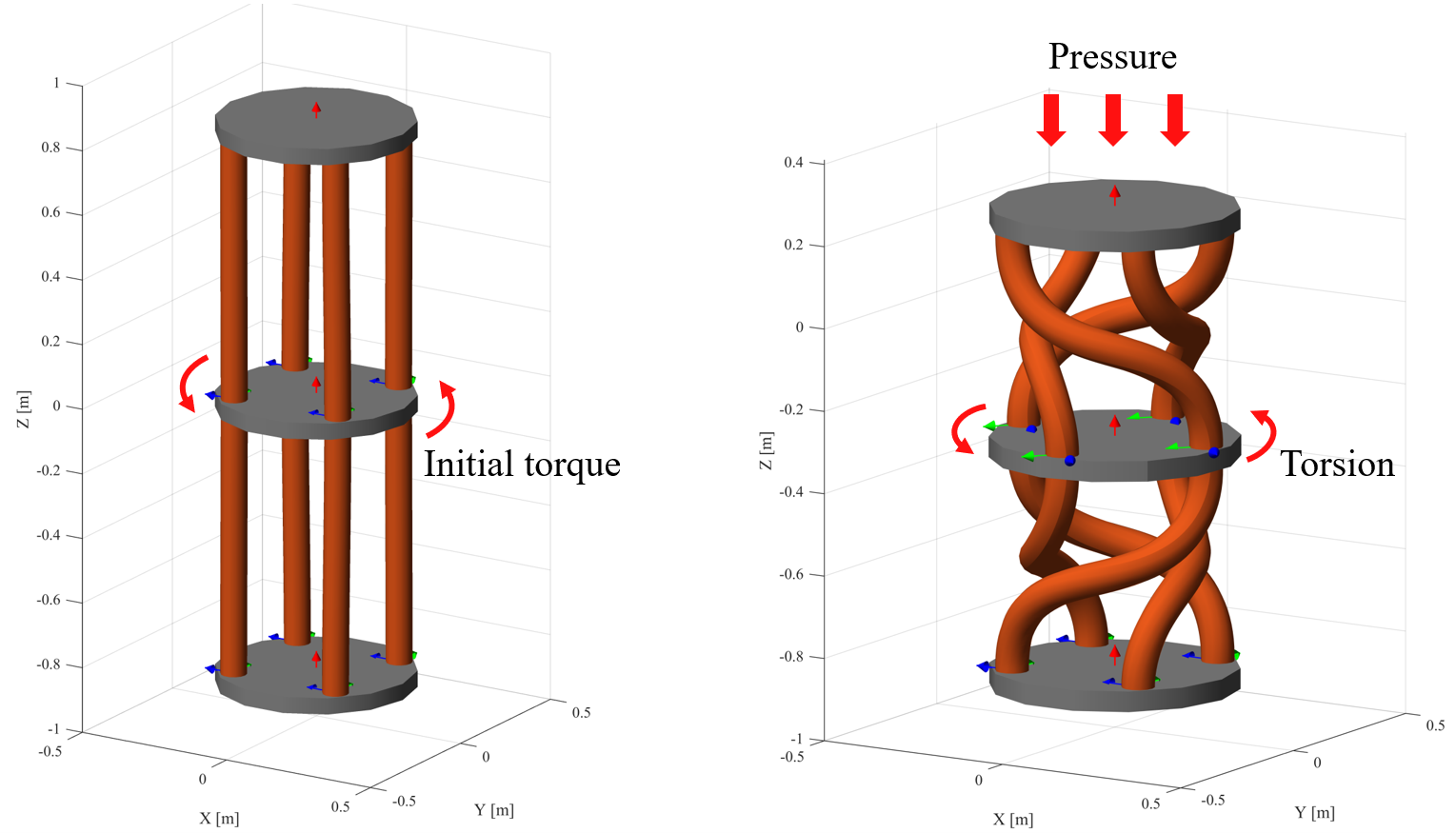}
	\caption{{Chiral structure under compressive loading.} The left image shows the initial unloaded state, in which straight vertical rods are symmetrically arranged between three horizontal circular plates. These rods are stress-free and maintain parallel alignment without deformation. The right image illustrates the deformed state under vertical compressive loading. The rods exhibit significant twist-buckling behavior, undergoing compression, bending, and torsion simultaneously.}
	\label{fig:twisting}
\end{figure}

We investigate a chiral parallel mechanism inspired by~\cite{fang2025large}, characterized by intrinsically coupled deformation modes—bending, twisting, and axial compression. As illustrated in~\Cref{fig:twisting}, the system comprises four compliant rods linking two rigid circular plates. Initially straight, vertically aligned, and free of residual stress, these rods exhibit a complex twist–buckling response under axial compression, enabling the structure to redistribute stress and dissipate energy.

To examine how different initial conditions influence the simulated behavior, we impose various levels of pre-twist torque on the upper plate prior to compression. The mechanism subsequently transitions through three representative configurations: the undeformed reference state, a mildly twisted pre-loaded state, and the final deformation shape under the applied compressive load.
\begin{table}[ht]
  \centering
  \small
  \caption{Physical Parameters}
  \label{tab:cosserat_params2}
  \begin{tabular}{l@{\hskip 115pt}c}
    \toprule
    \textbf{Parameter}        & \textbf{Value / Unit}            \\
    \midrule
    Length                   &  0.8 m                      \\
    Radius              & 4.0 cm                     \\
    Young’s modulus          &  $1.0\times10^{5}$ Pa     \\
    Possion's ratio        &  0.45\\
    \bottomrule
  \end{tabular}
\end{table}

In the reference configuration (\Cref{fig:twisting}, left), the mechanism comprises four straight, vertically aligned rods symmetrically linking a pair of circular plates. At this stage, the rods remain parallel and undamaged. A small vertical torque is subsequently applied to the upper disk to introduce a controlled pre-twist. By varying the magnitude of this applied torque, we generate a range of initial pre-twist conditions for comparison.

When the structure is compressed vertically (\Cref{fig:twisting}, right), the previously straight and unstressed rods deform into helical trajectories through twist-buckling, exhibiting a coupled response involving axial compression, bending, and torsion. This multi-mode deformation enables more uniform stress redistribution and substantially increases the structure’s capacity for elastic energy storage. Numerical simulations indicate that, under loading, the chiral mechanism achieves greater specific energy and energy density than conventional structures governed solely by axial buckling.

The material properties of the soft rods, including elastic modulus $E_s$, undeformed length $L_0$, and radius $r$, are listed in \Cref{tab:cosserat_params2}.
%
% ----------
%
% -------
% \subsubsection{Simulation results and analysis}
The simulation results obtained using our method are consistent with those reported in~\cite{fang2025large}, which used the finite element method, validating that our framework accurately captures the characteristic deformation behavior of chiral structures under compressive loading.
\begin{figure}[t]
	\centering
\includegraphics[width=0.95\textwidth]{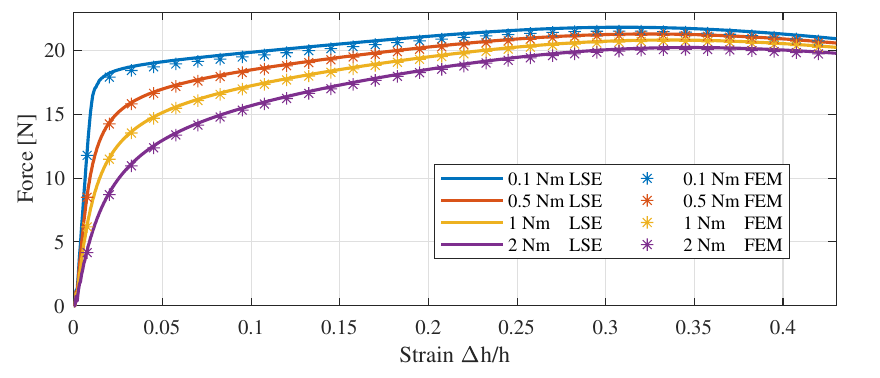}
	\caption{{Strain–force responses under different initial torques obtained with the proposed method using the linear strain element (LSE) and with the finite element method in COMSOL.} The plot illustrates the vertical compressive force versus strain relationship for chiral structures with varying initial torques (0.1~Nm, 0.5~Nm, 1~Nm, and 2~Nm). Lower initial torque results in a stiffer response with rapid force increase, while higher torque induces more pronounced twist-buckling and bending, allowing greater strain accommodation and enhanced energy absorption.}
	\label{fig:strain stress}
\end{figure}

\Cref{fig:strain stress} shows the simulated strain–force responses obtained with the proposed method using the linear strain element (LSE), alongside reference results from a finite element model implemented in COMSOL. The vertical axis in \Cref{fig:strain stress} corresponds to the magnitude of the reaction force at the base, while the horizontal axis represents the overall longitudinal compressive strain of the chiral structure, i.e., $\Delta h / h$, where $\Delta h$ denotes the vertical displacement of the upper rigid plate and $h$ is the initial height of the chiral assembly. As seen in the figure, the curves produced by our method are in excellent agreement with those of the COMSOL model.

\rv{In our setting, the pre-twist is intentionally small and mainly serves as a controlled perturbation that biases the structure towards the twist--buckling branch under subsequent vertical compression. As illustrated in \Cref{fig:strain stress}, once the vertical displacement is applied, the axial reaction rises rapidly at the beginning of the loading. In this stage, the response is dominated by axial compression: the axial strain and the corresponding axial resultant increase quickly, accounting for the main contribution to the reaction.}

\rv{As the compression proceeds, the axial reaction enters a plateau-like regime. Consistently, the axial strain and axial resultant tend to stabilize, indicating that further vertical shortening is no longer primarily accommodated by additional axial compression. Instead, geometric nonlinearity activates the coupled twist--bending mode, and the deformation progressively shifts to a torsion- and bending-dominated mechanism. In particular, the torsional and bending strain components increase and become the main channels to store the additional elastic energy, while the axial contribution remains comparatively stable. }

Overall, these results demonstrate that the proposed approach captures the transition between compression-dominated and twist--bending-dominated regimes without problem-specific parameter tuning, providing an efficient and reliable tool for simulating soft–rigid hybrid mechanisms and supporting the design and validation of multi–soft-rod structural systems.

\subsection{Gridshell}
% --------
%
Shell-like compliant structures are widely encountered in architectural structures and bio-inspired systems, including flexible protective skins, soft exoskeletons, bio-inspired wings, and related structures. Such systems combine low bending stiffness with distributed flexibility, requiring modeling approaches that remain effective under large geometric changes and non-trivial connectivity. To evaluate the proposed formulation in this context, we consider a fine rod lattice representation of a hemispherical shell.
\begin{figure}[!h]
    \centering
    \includegraphics[width=0.7\textwidth]{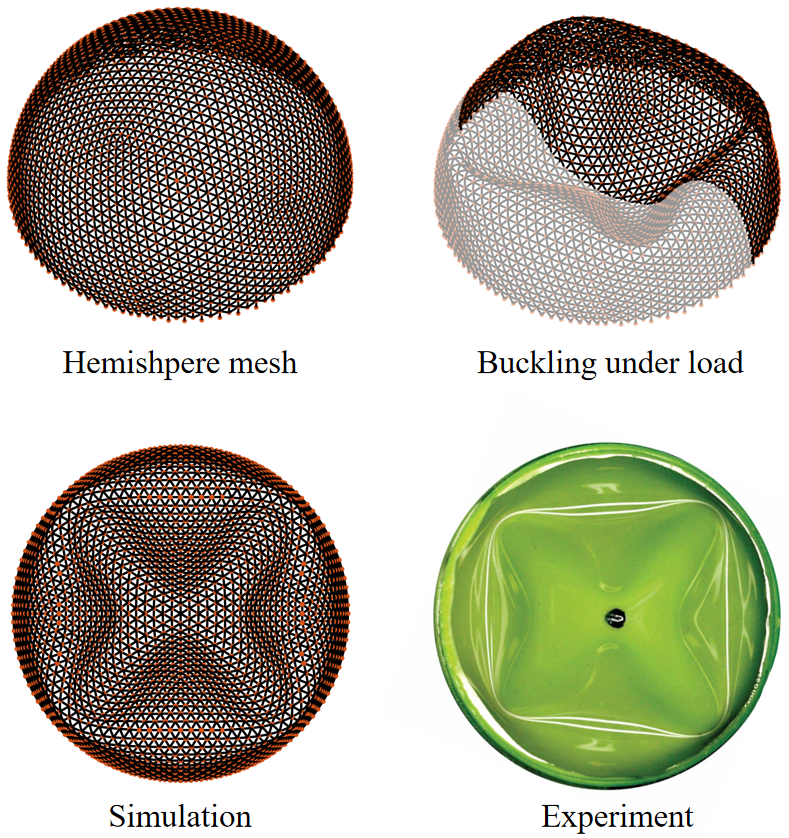}
    \caption{Hemispherical Cosserat-rod gridshell under a concentrated load. From left to right: initial mesh, intermediate deformation, post-buckling 
    configuration, and reference experimental image 
    (adapted from~\cite{nasto2013localization}).}
    \label{fig:hemisphere}
\end{figure}

As illustrated in \Cref{fig:hemisphere}, the configuration includes \(579\) nodes and \(1618\) Cosserat-rod elements, which form a triangular mesh over a hemisphere. The boundary at the equator is fixed, and a vertical load is applied at the pole. The components of the model and the solver configurations are the same as those utilized in the previous rod-network examples.

Under loading, the structure initially deforms in a manner consistent with shell-like global response and later develops a non-axisymmetric mode (\Cref{fig:hemisphere}), a behavior commonly observed in lightweight dome structures~\cite{nasto2013localization}. This transition emerges naturally from the equilibrium computation, without introducing additional constraints to the numerical algorithm.

Rod-lattice representations provide a discrete method for approximating the behavior of thin shells while maintaining a structured and sparse topology. As the lattice density increases, the mechanical response qualitatively approaches that of a continuous Cosserat-type shell model, providing a method for modeling shell-like soft structures with rod-based discretizations.

%
% ---------
\section{CONCLUSION}\label{sec.conc}
% ---------
%
This study presented a geometrically explicit formulation for Cosserat rods that combines Lie-group kinematics with local strain reconstruction. By treating nodal configurations on \( SE(3) \) as the primary degrees of freedom and recovering element-wise strains through a linear parameterization, the method merges the advantages of configuration-space and strain-based approaches within a unified framework. The resulting formulation preserves geometric consistency, avoids shear and membrane locking without ad hoc stabilization, and accommodates arbitrary rod networks and closed-loop arrangements through modular element assembly.

The numerical investigations demonstrate that the proposed approach captures large-deformation behavior with high accuracy, even with a relatively small number of elements. Its ability to scale seamlessly from single rods to complex lattices and gridshell-like structures highlights its relevance for the analysis of slender mechanisms, compliant structures, and soft robotic architectures.

\rv{Future developments will focus on extending the framework to dynamic simulations, incorporating reduced-order modeling for real-time applications, and integrating additional physical phenomena such as contact, actuation, and fluid--structure interactions. In particular, a key step toward dynamics is the consistent construction of inertial terms within the proposed element formulation: it requires computing the element kinetic energy and deriving its gradients with respect to the discrete degrees of freedom, in order to obtain the generalized inertial forces (and the associated Coriolis terms) in a form compatible with the \(SE(3)\) kinematics and the proposed pose--strain mapping. We will address this aspect in future work.} 

\rv{While this work focuses on rods, the proposed geometrically consistent formulation is not restricted to 1D and could be extended to 2D and 3D Cosserat models, e.g., Cosserat shells and micropolar continua. Investigating such generalizations, including their numerical robustness and potential locking behavior, is also an interesting direction for future work.}

\appendix
\section{Adjoint representation of the Lie algebra and Lie group}\label{notations}
	The adjoint representation of the Lie algebra $\mathfrak{se}(3)$ is given by	
	$$
	{\rm{ad}}_{{\xi}}= \left(\begin{matrix}
		\tilde{{\kappa}}&{0}\\\tilde{{\epsilon}}&\tilde{{\kappa}}
	\end{matrix}\right)\in \mathbb{R}^{6\times6}
	%, \ {\rm{ad}}_{\boldsymbol{\eta}}= \left(\begin{matrix}
		%	\tilde{\boldsymbol{w}}&\boldsymbol{0}_{3\times3}\\\tilde{\boldsymbol{v}}&\tilde{\boldsymbol{w}}
		%\end{matrix}\right)\in \mathbb{R}^{6\times6},
		$$
        where $\xi=[\kappa^\top \ \epsilon^\top]^\top\in \mathbb R^6$.
        
        The adjoint representation of the Lie group $SE(3)$ is given by	
		$${\rm{Ad}}_{{g}}= \left(\begin{matrix}
			{R}&{0}\\\tilde{{p}}{R}&{R}
		\end{matrix}\right)\in \mathbb{R}^{6\times6}.$$
        where $g\in SE(3)$.
% 
% ---------
\section{Tangent operator}
For $x\in SE(3)$, the tangent operator of the exponential map and its inverse are defined as follows:
\begin{equation}
	\mathrm{dexp}_x = \sum_{j=0}^{\infty} \frac{1}{(j+1)!} \, \mathrm{ad}^j_x, \quad 
	\mathrm{dexp}_x^{-1} = \sum_{j=0}^{\infty} \frac{B_j}{j!} \, \mathrm{ad}^j_x
\end{equation}
where \( B_j \) are the Bernoulli numbers. These series are typically truncated to achieve a specified level of accuracy. 
The first few Bernoulli numbers are given by \( B_0 = 1 \), \( B_1 = -\frac{1}{2} \), \( B_2 = \frac{1}{6} \), and \( B_3 = 0 \) \cite{hairer2006geometric}.
%
% --------
\section{Riemannian gradient on $SE(3)$}
\label{app:riem_grad}
This appendix summarizes the gradient computation on $SE(3)$ used in the static solver. We consider the left-multiplicative retraction
\begin{equation}
g(\epsilon)=g\,\exp(\epsilon\widehat{\zeta}),\qquad \zeta\in\mathbb{R}^6,
\label{eq:ret_se3}
\end{equation}
where $(\widehat{\cdot})$ and $(\cdot)^\vee$ denote the usual hat/vee operators on $\mathfrak{se}(3)$. Using a left-invariant metric, the (left-trivialized) Riemannian gradient $\mathrm{grad}^L_g U\in\mathbb{R}^6$ is defined by
\begin{equation}
\mathrm{D}U(g)\big[g\,\widehat{\zeta}\big]
\;=\;
\zeta^\top\,\mathrm{grad}^L_g U,
\qquad \forall\,\zeta\in\mathbb{R}^6.
\label{eq:left_grad_def}
\end{equation}

In this work, the element energy is written as an integral of a strain-energy density depending on the reconstructed linear strain field $\bar{\xi}(s)$, i.e.,
\begin{equation}
\begin{aligned}
U_e \;&=\;\frac{1}{2}\int_0^{h} (\xi-\xi_0)^\top K(\xi-\xi_0)\,ds \\
&= \; \frac{1}{2}h(\bar{\xi}-\xi_{0})^{\!\top}{K}(\bar{\xi}-\xi_{0})
    + \frac{1}{24}h^{3}\,\beta^{\!\top}{K}\beta.
\end{aligned}
\label{eq:strain_energy_density}
\end{equation}
with $K$ the Cosserat stiffness matrix in the body frame. $\bar{\xi}$ denote the average strain. Denote by $(g_a,g_b)\in SE(3)^2$ the two nodal poses of an element and by $\beta\in\mathbb{R}^6$ the strain-slope parameters of the LSE. We define the first-order variation of the  reconstructed linear strain as
\begin{equation}
\delta \bar{\xi}
\;=\;
J_1\,\delta \zeta_a \;+\; J_2\,\delta \zeta_b \;+\; J_3\,\delta \beta,
\label{eq:delta_xi_chain}
\end{equation}
where $\delta \zeta_a,\delta \zeta_b\in\mathbb{R}^6$ are defined by $\delta g_a=g_a\widehat{\delta\zeta_a}$ and $\delta g_b=g_b\widehat{\delta\zeta_b}$ (consistent with \eqref{eq:ret_se3}), and $J_1,J_2,J_3$ are the Jacobians given in the main text.

Combining \eqref{eq:left_grad_def} and \eqref{eq:strain_energy_density} yields
\begin{equation}
\delta U_e
=\delta \bar{\xi}^\top{K}(\bar{\xi}-\xi_{0})
    + \frac{1}{12}h^{3}\,\delta {\beta}^\top{K}\beta
=
\delta \zeta_a^\top R_{1} \;+\; \delta \zeta_b^\top R_{2} \;+\; \delta \beta^\top R_{3},
\label{eq:deltaU_to_residual}
\end{equation}
Taking \eqref{eq:delta_xi_chain} into \eqref{eq:deltaU_to_residual}, one can get the formulation of the Riemannian gradients:
\begin{equation}
\begin{aligned}
R_1 &= h\,J_1^{\!\top}\, K\,(\bar\xi-\xi_0),
\\[2pt]
R_2 &= h\,J_2^{\!\top}\, K\,(\bar\xi-\xi_0),
\\[2pt]
R_3 &= h\,J_3^{\!\top}\, K\,(\bar\xi-\xi_0)
      + \tfrac{1}{12}h^3\, K\,\beta.
\end{aligned}
\end{equation}

\section{Derivation of Eq. \eqref{eq:magnus4}}
\label{app:eq5}

On one element $s\in[0,h]$, the right-multiplicative kinematics reads
\begin{equation}
g'(s)=g(s)\,\widehat{\xi}(s), \qquad g(0)=g_a,\quad g(h)=g_b .
\end{equation}
It follows that $g_a^{-1}g_b=\exp(\widehat{\Omega})$ for some $\Omega\in\mathbb{R}^6$.
We approximate $\widehat{\Omega}$ by a fourth-order Magnus truncation written in the compact
moment form \cite{blanes2006fourth}:
\begin{equation}
\widehat{\Omega}\;\approx\;\int_0^h\widehat{\xi}(s)\,ds
\;-\;\Bigg[\frac{1}{h}\int_0^h\Big(s-\frac{h}{2}\Big)\widehat{\xi}(s)\,ds,\ \int_0^h\widehat{\xi}(s)\,ds\Bigg].
\end{equation}
Under the linear strain assumption $\xi(s)=\bar\xi+(s-\tfrac{h}{2})\beta$, direct integration gives
$$\int_0^h\widehat{\xi}(s)\,ds=h\,\widehat{\bar\xi}, \ \ \frac{1}{h}\int_0^h(s-\frac{h}{2})\widehat{\xi}(s)\,ds=\frac{h^2}{12}\widehat{\beta},$$
hence
\begin{equation}
{\Omega}
= h\,{\bar\xi}-\frac{h^3}{12}\big[{\beta},{\bar\xi}\big]
\;=\;
{\Big(hI-\frac{h^3}{12}\mathrm{ad}_{\beta}\Big)\bar\xi},
\end{equation}
which yields Eq.~\eqref{eq:magnus4} in the main text.

\bibliographystyle{elsarticle-num}
\bibliography{biblio.bib}

% %% For numbered reference style

%   \textit{\LaTeX: a document preparation system},
%   Addison Wesley, Massachusetts,

% \end{thebibliography}
\end{document}